\documentclass[aos,MSNbibl,nameyear,dvips]{arximspdf}
\usepackage{graphicx}
\doi{10.1214/14-AOS1239} 
\volume{42}
\issue{5}
\pubyear{2014}
\firstpage{1850}
\lastpage{1878}
\docsubty{FLA}

\makeatletter
\newcommand{\widebar}{\overline}
\newcommand{\TV}{\mathrm{TV}}
\renewcommand{\mid}{|}
\newtheorem{theorem}{Theorem}
\newtheorem{lemma}{Lemma}
\newproclaim{remark}{Remark}
\newtheorem{corollary}{Corollary}
\newtheorem{proposition}{Proposition}
\newproclaim{example}{Example}
\newproclaim{ass}{Assumption}
\makeatother

\begin{document}
\begin{frontmatter}

\title{The Bernstein--von Mises theorem and
nonregular~models\thanksref{T1}}
\runtitle{Bernstein--von Mises theorem and nonregular models}

\begin{aug}
\author[a]{\fnms{Natalia A.}~\snm{Bochkina}\corref{}\ead[label=e1]{N.Bochkina@ed.ac.uk}}
\and
\author[b]{\fnms{Peter J.}~\snm{Green}\ead[label=e2]{P.J.Green@bristol.ac.uk}}
\runauthor{N. A. Bochkina and P. J. Green}
\affiliation{University of Edinburgh and Maxwell Institute,\\
and
University of Technology, Sydney and University of Bristol}
\address[a]{School of Mathematics\\
University of Edinburgh\\
Edinburgh EH9 3JZ\\
United Kingdom\\
\printead{e1}}
\address[b]{School of Mathematics\\
University of Bristol\\
Bristol BS8 1TW\\
United Kingdom\\
\printead{e2}}
\end{aug}
\thankstext{T1}{Both authors acknowledge financial
support for research visits provided by the EPSRC-funded SuSTaIn
programme at Bristol University.}

\received{\smonth{7} \syear{2013}}
\revised{\smonth{5} \syear{2014}}

%
\begin{abstract}
We study the asymptotic behaviour of the posterior distribution in a
broad class of statistical
models where the ``true'' solution occurs on the boundary of the
parameter space. We show that
in this case Bayesian inference is consistent, and that the posterior
distribution has not only
Gaussian components as in the case of regular models (the
Bernstein--von Mises theorem) but also
has Gamma distribution components whose form depends on the behaviour
of the prior distribution
near the boundary and have a faster rate of convergence. We also
demonstrate a remarkable property of Bayesian inference, that for some models,
there appears to be no bound on efficiency of estimating the unknown
parameter if it is on the boundary of the parameter space.
We illustrate the results on a problem from emission tomography.
\end{abstract}

%
\end{frontmatter}

\section{Introduction}

The asymptotic behaviour of Bayesian methods has been a long-standing
topic of interest, including approximation of the posterior
distribution and questions that are important from a frequentist point
of view, such as consistency, efficiency and coverage of Bayesian
credible regions.
For instance, for correctly specified regular finite-dimensional models
with $n$ independent observations, these properties are captured by the
Bernstein--von Mises theorem that implies that the posterior
distribution can be approximated in a $1/\sqrt{n}$ neighbourhood of
the true value of the parameter by a Gaussian distribution with
variance determined by the Fisher information. More generally, the
Bernstein--von Mises theorem holds for dependent observations if the
likelihood satisfies local asymptotic normality (LAN) conditions
[\citet{LeCam,LeCamYang}]. A~total variation distance version of
the theorem was derived by \citet{vdVaart}.
This theorem implies that the prior has no asymptotic influence on the
posterior, that posterior inference is consistent and efficient in the
frequentist sense, and that posterior credible regions are
asymptotically the same as frequentist ones.

One of the key assumptions of the Bernstein--von Mises (BvM) theorem is
that the ``true'' value of the parameter is an interior point of the
parameter space.
However, for many problems, including our motivating example of a
Poisson inverse problem in tomography, and, more generally for the
class of models we consider, this assumption of the BvM theorem does
not hold. For the tomography example, the unknown parameter is a vector
of tracer concentrations, which are nonnegative and can be zero.

The situation where the unknown parameter can be on the boundary of the
parameter support has been addressed in the frequentist literature by
studying the asymptotic distribution of the maximum likelihood
estimator [\citet{SelfLiang,Moran}, among others]; however it has
been studied very little under the Bayesian approach. \citet
{HDudley} investigated the asymptotic behaviour of the posterior
probability of the unknown parameter belonging to a half-space
$\mathcal{H}$ for a regular correctly specified model, where they
found that if the true value of the parameter belongs to the complement
of $\mathcal{H}$, then the posterior probability of half-space
$\mathcal{H}$ goes to zero much faster, namely at least at rate $n$
rather than at the standard parametric rate $\sqrt{n}$ ($n$ here is
the sample size), and there is an exponential upper bound on this
posterior probability.
Also, \citet{Erkanli} gave a formula for calculating the
expectation of a smooth functional of a 3-dimensional posterior
distribution where the unknown parameter is on a smooth boundary.

In this paper, we extend the Bernstein--von Mises theorem by relaxing
the assumption that the ``true'' value of the parameter is interior to
the parameter space, in a finite-dimensional setting. We consider a
broad class of probability distributions for the data and allow the
prior distribution to be improper and to have zero or infinite density
on the boundary. A key model assumption is that the ``true'' value of
the parameter minimises a generalised Kullback--Leibler distance. There
is no assumption of any finite moments. We will show that for these
models the consequences of relaxing this assumption are twofold:
firstly, the convergence is faster, at least at rate $n$, if the
``true'' parameter is on the boundary, and secondly, the limit of the
posterior distribution has non-Gaussian components.

We motivate our study by presenting in Section~\ref{secMotivation} an
inverse problem from medical imaging; Section~\ref{secSetup}
establishes the class of models we study.
In Section~\ref{secMainResult} we state the result on the local
behaviour of the posterior distribution in a neighbourhood of the limit
that is formulated as a modified Bernstein--von Mises theorem, discuss
the assumptions, and give a nonasymptotic version of the result. In
Section~\ref{secBvMexamples} we illustrate the application of our
analogue of the BvM theorem for various examples including the problem
of variance estimation in mixed effects models, and discuss the choice
of prior distribution.
We discuss issues in using the approximation of the posterior
distribution in practice and apply it to data from the motivating
example in Section~\ref{secSPECTexample}. We conclude with a
discussion. All proofs are deferred to the \hyperref[append]{Appendix}.

\section{Motivating example}\label{secMotivation}

\subsection{Single photon emission computed tomography}\label{secspect}

Single photon emission computed tomography (SPECT) is a medical
imaging technique in which a radioactively-labelled tracer, known to
concentrate in the tissue to be imaged, is introduced into the subject. Emitted
particles are detected in a device called a gamma camera,
forming an array of counts. Tomographic
reconstruction is the process of inferring the spatial pattern of
concentration of the tracer in the tissue from these
counts. The Poisson linear model
%
\begin{equation}
\label{eqPoisError} \mathcal{T}Y_i\mid\theta\sim\operatorname{Poisson}(
\mathcal{T}A_i \theta),\qquad i=1,\ldots, n,\mbox{ independently},
\end{equation}
comes close to reality for the SPECT problem (there are some
dead-time effects and other artifacts in recording). Here
$\theta=\{\theta_j\},j=1,2,\ldots,p$
represents the spatial distribution of the tracer,
typically discretised on a grid, with $\theta_j \geq0$ for all $j$,
$Y=\{Y_i\}$
the array of the rate of detected photons per time unit, also
discretised by the
recording process, and $\mathcal{T}$ is the exposure time for photon
detection. The
$n \times p$
array $A=(A_{ij})$ with rows $A_i$ quantifies the emission,
transmission, attenuation, decay, and recording process; $A_{ij}$
is the mean number of photons recorded at $i$ per unit
concentration per unit time at pixel/voxel $j$, and is nonnegative.
In some methods of reconstruction, elements of the matrix $A$ are
modelled as discretised values of the Radon transform.

Since Poisson distributions form an exponential family, this model can
be seen as a generalised linear model
[\citet{NelderW}], with identity link function and dispersion
$1/\mathcal{T}$; see also Example~\ref{exPoissonM} in Section~\ref{secExamplesBoundary}.

We formalise the notion of small-noise limit for this Poisson
model in a practically-relevant way, by supposing that the
exposure time for photon detection becomes large, that is, letting
$\mathcal{T}\to\infty$.

The ``true image''
$\theta^\star$ in emission tomography corresponds to a physical reality,
the discretised spatial distribution of concentration of the tracer.
Since this is nonnegative, we impose the constraints $\theta\in
\Theta= [0,\infty)^p\subset\mathbb{R}^p$.

Unless $p$ is too large, that is, the spatial resolution of $\theta$
is too fine, the matrix $A$ is normally of full rank $p$, and hence the
inverse problem is well posed (although it may be ill-conditioned); see
\citet{JSPET} for eigenvalues of the Radon transform.

See \citet{Green} for further detail about this model, and an
approach based on EM estimation for MAP
reconstruction of $\theta$, in a Bayesian formulation in which spatial
smoothness of the solution is promoted by using a pairwise
difference Markov random field prior.

\subsection{Prior distribution}\label{secpriors}

From the beginning of Bayesian image analysis
[\citet{GemanG,Besag}], use has been made of
Markov random fields as prior distributions for image scenes that
express generic,
qualitative beliefs about smoothness, yet do not rule out abrupt
changes for real discontinuities (e.g., at
tissue type boundaries in the case of medical imaging).

The prior distribution we consider for the SPECT model is a log cosh
pairwise-interaction Markov
random field [\citet{Green}],
%
\begin{equation}
\label{eqdefLogCosh} p(\theta) \propto\exp \biggl(-\frac{\zeta(1+\zeta)}{2\gamma^2} \sum
_{j\sim j'} \log\cosh \biggl(\frac{\theta_j-\theta_{j'}}{\zeta} \biggr) \biggr),
\qquad\theta\in\Theta,
\end{equation}
where $j\sim j'$ stands for $j$ and $j'$ being
neighbouring pixels. In this paper the parameters $\zeta$ and $\gamma
$ are considered to be fixed.

This model has some attractive properties. While giving less penalty to
large abrupt changes in $\theta$ compared to
the Gaussian, it remains log-concave. It bridges the extremes
$\zeta\to\infty$, the Gaussian pairwise-interaction prior, and
$\zeta=0$,
the corresponding Laplace pairwise-interaction model, sometimes called
the ``median prior.''

This distribution is improper since it is invariant to perturbing
$\theta$ by an arbitrary
additive constant, but leads to a proper posterior distribution as long
as $\sum_{j} A_{ij} \neq0$ for some $i$.

\subsection{Nonstandard features of the SPECT model}

The Bayesian model for SPECT has three nonstandard features: (a) the
true image $\theta^\star$ can lie on the boundary of the parameter
space $[0,\infty)^p$; (b) if $A_i \theta^\star=0$ for some $i$, then
the distribution of the corresponding $Y_i$ degenerates to a point mass
at 0; (c) the prior distribution is not proper.

In the next section we formulate a model that includes the Bayesian
SPECT model as a particular case. The approximate behaviour of the
posterior distribution of $\theta$ for large $\mathcal{T}$ is
investigated in Section~\ref{secSPECTexample}.

\section{Model formulation}\label{secSetup} 

\subsection{Likelihood}\label{secLik}

We now list assumptions on the distribution of the observable responses
$Y$, taking values in $\mathcal{Y}\subseteq\mathbb{R}^n$; it has
density (with respect to Lebesgue or counting measure) denoted by
$p_{\sigma} ( y \vert \theta)$ for $\theta\in\Theta\subset
\mathbb{R}^p$.\vadjust{\goodbreak}
These assumptions are expressed in terms of the scaled log-likelihood
defined by
\[
\ell_{y,\sigma}(\theta) = \sigma^2\log p_{\sigma} ( y
\vert \theta). %
\]
As we shall see from the assumptions, $\sigma$ is related to the level
of noise, and we are interested in the case where $\sigma$ is small.
We assume that the ``true'' value of the unknown parameter that
generated the data is $\theta^\star\in\Theta$, and denote the true
probability measure of $Y$ by $\mathbb{P}_{\theta^\star, \sigma}$.
Below, where it does not lead to ambiguity, we will omit the index
$\sigma$ to simplify the notation and will write $\ell_{y}(\theta)$
and $\mathbb{P}_{\theta^\star}$.

\renewcommand{\theass}{M}
\begin{ass}\label{assM}
(1) For $Y \sim p_{\sigma} ( y \vert \theta^\star)$,
there exists a deterministic function $\ell^\star(\theta)\dvtx  \Theta
\to\mathbb{R}$ such that for all $\theta\in\Theta$,
\[
\forall\varepsilon>0 \qquad\mathbb{P}_{\theta^\star, \sigma} \bigl( \bigl|
\ell_{Y(\omega)}(\theta) - \ell^\star(\theta)\bigr|>\varepsilon \bigr) \to0
\qquad\mbox{as } \sigma\to0. 
\]

(2) The function $\ell^\star(\theta)$ has a unique maximum
over $\Theta$ at $\theta= \theta^\star$.
\end{ass}
Further assumptions on $\ell_{y,\sigma}(\theta)$ are given in
Section~\ref{secQuadraticApprox}.

Assumption~\ref{assM} is satisfied for a wide class of models, in particular for
models with independent identically distributed (i.i.d.) observations
with $\sigma^2 =1/n$ and for distributions from the exponential family
in canonical form with dispersion $\sigma^2\to0$, that are discussed below.

The function $\ell^\star(\theta)$, defined in Assumption~\ref{assM}(1), can be
viewed as the limit of the negative Kullback--Leibler (${\mathcal
{KL}}$) distance, rescaled by $\sigma^2$, between distributions with
densities $p (\cdot\mid\theta)$ and $p (\cdot\mid\theta^\star
)$, that was used, for instance, in \citet{PetroneEB} and
\citet{Barron99}.
For i.i.d. models, $\ell^\star(\theta)$ is the negative
Kullback--Leibler distance based on a single observation, and for
generalised linear models
$\ell^\star(\theta)$ is the log-likelihood for ``noise-free'' data.
Assumption~\ref{assM}(2) states that this generalised Kullback--Leibler distance
is minimised at the ``true'' value $\theta^\star$, as holds for the
usual ${\mathcal KL}$ distance.
Assumption~\ref{assM}(2) has been used by other authors, for instance, in the
context of hidden Markov models by \citet{DoucHMM} where it was
called the \textit{identifiability assumption}, and a finite sample
analogue of this assumption was used in the context of a misspecified
model by \citet{Spokoiny2012}.
This assumption holds for some models where the parameter set $\Theta$
is not open and thus where the true value of the parameter $\theta
^\star$ can be on the boundary of $\Theta$; see Example~\ref{exPoissonM}.
These assumptions are satisfied for the tomography model discussed in
Section~\ref{secMotivation} where the unknown tracer image $\theta
^\star$ can have zero intensity values in some pixels, as shown in
Section~\ref{secExamplesBoundary}.

Next we show that Assumption~\ref{assM} is satisfied for two important classes
of models, generalised linear models, and i.i.d. models, including the
case when $\theta^\star$ is on the boundary of $\Theta$.

\subsection{Generalised linear models}\label{secExamplesBoundary}

In the generalised linear models of \citet{NelderW}, an important
class of nonlinear
statistical regression problems, responses $y_i$, $i=1,2,\ldots,n$ are drawn
independently from a one-parameter exponential family of distributions
in canonical form, with density or probability function
\[
p_{\sigma}(y\mid\eta) = \exp \Biggl(\sum_{i=1}^n
\biggl[\frac
{y_ib(\eta_i)-c(\eta_i)}{\sigma^2} +d(y_i,\sigma) \biggr] \Biggr), %
\]
using the mean parameterisation, for appropriate functions
$b$, $c$, and $d$ characterising the particular distribution
family. The parameter $\sigma^2$ is a common dispersion parameter
shared by all responses. Assuming that functions $b(\cdot)$ and
$c(\cdot)$ are twice differentiable, the expectation of this distribution
is $\mathbb{E}(Y_i)=\eta_i=c'(\eta_i)/b'(\eta_i)$, and the variance
is $\operatorname{Var}(Y_i)= \sigma^2 [ c''(\eta_i)b'(\eta
_i)-c'(\eta_i)
b''(\eta_i)] /[b'(\eta_i)]^3$.
This implies that the random variable $Y$ converges in probability to a
finite deterministic limit $y^\star= \mathbb{E}Y$ as $\sigma\to0$
and that the dispersion $\sigma^2$ is related to the noise level of
the observations.

Firstly consider the case $\theta= \eta$. Then, $\ell_{Y}(\theta)$
is linear in $Y$, and hence it converges to $\ell_{y^\star}(\theta)$
in probability as $\sigma\to0$. Therefore, Assumption~\ref{assM}(1) is
satisfied with $\ell^\star(\theta) = \ell_{y^\star}(\theta)$.
If $\nabla\ell^\star(\theta^\star)=0$ and the Hessian,
which is diagonal, has negative entries, then
$\theta^\star$ uniquely maximises $\ell^\star(\theta)$; that is,
Assumption~\ref{assM}(2) is satisfied. If $\theta^\star$ is on the boundary and
the gradient is nonzero, see Examples~\ref{exPoissonM} and~\ref
{exBinomM} below.

Now consider a generalised linear model with $\eta= A\theta$ and
matrix $A$ such that $A^T A$ is of full rank, that is, such that the
likelihood is identifiable with respect to parameter $\theta$. In this
case, Assumption~\ref{assM} holds with $\theta^\star= (A^T A)^{-1} A^T y^\star$.
The tomography example given in Section~\ref{secMotivation} belongs
to this class of models, with $\sigma^2 = \mathcal{T}^{-1}$, $b(\eta
_i) = \log\eta_i$, $c(\eta_i) = \eta_i$, 
and $\Theta= [0,\infty)^p$.

Now we show that Assumption~\ref{assM}(2) is satisfied when $\theta^\star$ is
on the boundary of~$\Theta$ for some distributions from the
exponential family.
%
\begin{example}\label{exPoissonM}
Consider the Poisson distribution $Y/\sigma^2 \sim\operatorname
{Poisson}(\eta/\sigma^2)$ with $\eta\geq0$. The scaled
log-likelihood for $\eta$ is
$\ell_y(\eta) = y\log\eta- \eta$.
If $Y$ is generated with $\eta=0$, then we observe $y=0$ with
probability 1, so in this case the scaled log-likelihood for $\eta$ is
always $-\eta$, which is maximised over $\eta\geq0$ at $\eta=0$,
that is, the true value of $\eta$.
\end{example}

\begin{example}\label{exBinomM}
For the Binomial distribution $Y\sim\operatorname{Bin}(n, \eta)$,
the scaled log-likelihood for $\eta\in[0,1]$ is
$
\ell_y(\eta) = [y\log ( \eta ) + (n-y)\log(1-\eta)]/n$.
If the true value of $\eta$ is $1$, then $\mathbb{P}(Y=n)=1$ and the
scaled log-likelihood for $\eta$ is $\ell_y(\eta) = \log( \eta)$,
which is maximised over $[0,1]$ at $\eta=1$, so that again we recover
the true value, and Assumption~\ref{assM}(2) is satisfied for this model.

The same holds for the other boundary point $\eta=0$, and also for
multinomial and negative binomial distributions.
\end{example}

\subsection{I.I.D. models}

Let $Y_1,\ldots, Y_n$ be independent identically distributed random
variables where the density or probability mass function of $Y_i$ is
$p(y_i \mid\theta) = C_{y_i} \exp\{\ell_{y_i}(\theta)\} $, with
unknown parameter $\theta\in\Theta\subset\mathbb{R}^p$ where $p$
is finite and independent of $n$.
Here, $\sigma^2 = 1/n$ and $\ell_y(\theta) = n^{-1} \sum_{i=1}^n
\ell_{y_i}(\theta)$.
In this case, $\ell_{Y_i}(\theta)$ are i.i.d. random variables, so,
as $n\to\infty$, Assumption~\ref{assM}(1) is satisfied under the conditions of
the weak law of large numbers for the random variable $\ell
_{Y_i}(\theta)$, for all~$\theta$, which implies that there exists
$\ell^\star(\theta)$ such that
$\ell_Y(\theta)$ converges in probability to $\ell^\star(\theta)$
as $n\to\infty$.
If $\mathbb{E}[\ell_{Y}(\theta)]$ exists for all $\theta\in\Theta
$, then $\ell^\star(\theta) = \mathbb{E}[\ell_{Y_i}(\theta)]$,
equal to the negative Kullback--Leibler distance between the
distributions with densities $p(\cdot\mid\theta)$ and $p(\cdot\mid
\theta^\star)$, and then Assumption~\ref{assM}(2) holds.
For instance, it is easy to check that Assumption~\ref{assM} is satisfied for
i.i.d. Cauchy random variables $Y_i$ with $\ell_{Y_i}(\theta) = \log
 (1+ (Y_i - \theta)^2  )$ and $\theta\in\Theta\subseteq
\mathbb{R}$.

\subsection{Bayesian formulation}

We adopt a Bayesian paradigm, using a $\sigma$-finite prior measure
$\pi(d\theta)$ on $\Theta$.
Thus the posterior distribution satisfies
%
\begin{equation}
\label{eqPosterior} \pi(d \theta\vert y) \propto\exp\bigl( \ell_y(
\theta)/\sigma^2 \bigr) \pi(d\theta), \qquad\theta\in\Theta.
\end{equation}
Here we do not assume that the prior distribution is proper, nor do we
assume that its density is bounded away from 0 and infinity on the
boundary of $\Theta$; see Assumption~\ref{assP} in Section~\ref{secQuadraticApprox}.

\section{The analogue of the Bernstein--von Mises theorem}\label{secMainResult}

\subsection{Notation and assumptions}\label{secQuadraticApprox}

We shall use the default norms $\|z\| = \|z\|_2$ for both vectors and
matrices. If the appropriate derivatives exists, define the gradient
$\nabla f(\theta)$ of a function $f$ on $\Theta$ as a vector of
partial derivatives (one-sided if $\theta$ is on the boundary of
$\Theta$), and $\nabla^2 f(\theta)$ is a matrix of second
derivatives of $f$ (again, one-sided if $\theta$ is on the boundary of
$\Theta$). We use notation $\theta_{S}$ to define the vector \mbox{$(\theta
_{j}, j\in S)$} for $S\subset\{1,\ldots, p\}$, a convention which
also applies to the gradient $\nabla$, that is, $\nabla_{S} f(\theta
)= (\nabla_{j} f(\theta), j\in S)$. We denote a submatrix $\Sigma
$ indexed by subsets $S, J$ by $\Sigma_{S, J} = (\Sigma_{ij}, i\in
S, j\in J)$; this also applies to the matrix of second derivatives,
so we can write $\nabla^2_{S,J} f(\theta)$ to denote the
corresponding submatrix.

We use $A\mathcal{X}+ x_0= \{Ax+x_0, x\in\mathcal{X}\}$ to denote
the image of an affine transformation of the set $\mathcal{X}$ given
matrix $A$ and vector $x_0$.

The limit of the posterior distribution has a different character in
different directions, and we need to partition the index set $\{
1,2,\ldots,p\}$ of $\theta$ accordingly. Let
\[
S_0=\bigl\{j\dvtx  \nabla_j \ell^\star\bigl(
\theta^\star\bigr) = 0\bigr\}\quad\mbox {and}\quad S_1=\bigl
\{j\dvtx  \nabla_j \ell^\star\bigl(\theta^\star\bigr)
\neq0\bigr\}, %
\]
with dimensions $p_0$ and $p_1=p-p_0$, respectively.
We partition $S_0$ further:
\begin{eqnarray*}
S_0^\star&=&\bigl\{j\dvtx  \nabla_j
\ell^\star\bigl(\theta^\star\bigr) = 0\mbox{ and }
\theta^\star_j = 0\bigr\} \quad\mbox{and}
\\
S_0\setminus S_0^\star&=& \bigl\{j\dvtx
\nabla_j \ell^\star\bigl(\theta^\star\bigr) = 0
\mbox{ and } \theta^\star_j \neq0\bigr\},
\end{eqnarray*}
with dimensions $p_0^\star$ and $p_0-p_0^\star$ {where $\theta^\star
_j = 0$ corresponds to $\theta^\star$ being on the boundary of
$\Theta$; see Assumption~\ref{assB}(1) below}.

We then introduce a permutation of coordinates of $\theta$, defined by
any matrix $U$ that maps $S_0 \setminus S_0^\star$ to the first
$(p_0-p^\star_0)$ coordinates, $S_0^\star$ to the next $p^\star_0$,
and $S_1$ to the last $p_1$. The first $p_0$ rows of $U$ will be
denoted $U_0$ and the remainder $U_1$. We denote the index set $\{
p_0-p_0^\star+1,\ldots,p_0\}$ by $T_0^\star$ which is the image of
$S_0^\star$ under the map defined by $U$. Note that $\theta^\star
_j=0$ for all $j\in S_0^\star\cup S_1$ (for $j\in S_1$, this is given
by Lemma~\ref{lemopt} below), so this set describes the coordinates of $\theta
^\star$ that lie on the boundary; in the case of $S_0^\star$ the
gradient is also zero in this direction.

We introduce the notation $\mathcal{PTN}_{p_0}(a_0,\Omega_{00}^{-1},
p_0^\star,\alpha_0)$ for a polynomially-tilted multivariate Gaussian
distribution truncated to $\mathcal{V}_0 = \mathbb{R}^{p_0-p^\star
_0}\times\mathbb{R}_+^{p^\star_0}$, for which the corresponding
measure of any measurable $\mathcal{B}\subset\mathcal{V}_0$ is
defined by
%
\begin{eqnarray}\label{eqModifGauss}
&& \mathcal{PTN}_{p_0}\bigl(\mathcal{B};a_0,
\Omega_{00}^{-1}, p_0^\star,
\alpha_0\bigr)
\nonumber\\[-1pt]\\[-15pt]
&&\qquad  = \frac{\int_{\mathcal{B}} \prod_{j\in T_0^\star}
x_j^{\alpha_{0,j} -1}
e^{ - (x - a_0 )^T \Omega_{00}(x - a_0 ) /2 } \,dx}{
\int_{\mathcal{V}_0} \prod_{j\in T_0^\star} x_j^{\alpha_{0,j} -1}
e^{ - (x - a_0 )^T \Omega_{00}(x - a_0) /2 } \,dx},\nonumber
\end{eqnarray}
where $a_0 \in\mathbb{R}^{p_0}$, $\Omega_{00}$ is a $p_0\times p_0$
positive definite matrix, and $\alpha_0=(\alpha_{0,j})_{j\in
T_0^\star} \in(0,\infty)^{p^\star_0}$. $\alpha_0$ could also be
interpreted as a $p_0$-dimensional vector whose first $p_0-p_0^\star$
coordinates are irrelevant.
Note that this distribution is Gaussian if $p^\star_0=0$, and
truncated Gaussian if $p^\star_0\neq0$ and $\alpha_{0,j}=1$ for all $j$.

For $\alpha, a>0$, $\Gamma(\alpha, a)$ denotes the Gamma
distribution with density $p(x) = a^{\alpha} x^{\alpha-1}
e^{-ax}/\Gamma(\alpha)$, $x>0$,
and $\Gamma(dx; \alpha, a)$ the corresponding probability measure.

In addition to Assumption~\ref{assM} (Section~\ref{secLik}), we make the
following assumptions.
They make use of the following neighbourhoods of $\theta^\star$:
%
\begin{equation}
\label{eqDefBdelta} \Theta^\star(\delta) = \bigl\{\theta\in\Theta\dvtx  U\bigl(
\theta-\theta ^\star\bigr) \in B_{2, p_0}(0, \delta_0)
\times B_{\infty, p_1}(0, \delta _1)\bigr\},
\end{equation}
where $\delta=(\delta_0, \delta_1)$, $\delta_0, \delta_1 >0$ and
$B_{q,s}(z_0, r) = \{z\in\mathbb{R}^s\dvtx  \|z-z_0\|_q < r \}$.

\renewcommand{\theass}{B}
\begin{ass}[(On boundary of $\Theta$, $\partial \Theta$)]\label{assB}
(1) $\Theta\subseteq[0,\infty)^p$ and $\Theta\cap\partial\Theta
\subseteq\bigcup_{j=1}^p \{\theta\in\Theta\dvtx  \theta_j =0\}$.

(2) $U(\Theta-\theta^\star) \supseteq(-c_0,c_0)^{p_0-p_0^\star}
\times[0,c_0)^{p_0^\star} \times[0,c_1)^{p_1}$ for some $c_0,c_1>0$.
\end{ass}

\renewcommand{\theass}{S}
\begin{ass}[(Smoothness in $\theta$)]\label{assS}
There exist $\delta_0, \delta_1 >0$ depending on $\sigma$ such that:
\begin{longlist}[(2)]
\item[(1)] $\delta_0 \to0$, $\delta_1 \to0$, $\delta_0/\sigma\to\infty
$, $\delta_1/\sigma^2 \to\infty$ as $\sigma\to0$.

\item[(2)] For all $\theta\in\Theta^\star(\delta)$, $\nabla\ell^\star
(\theta)$, $\nabla\ell_Y(\theta)$ and $\nabla^2_{S_0, S_0} \ell_Y
(\theta)$ exist $\mathbb{P}_{\theta^\star, \sigma}$-almost
everywhere, for small enough $\sigma$.

\item[(3)] For any $\varepsilon>0$,
\[
\mathbb{P}_{\theta^\star, \sigma} \Bigl( \sup_{\theta\in\Theta
^\star(\delta)} \bigl\|\nabla
\ell_{Y(\omega)}(\theta) - \nabla\ell ^\star\bigl(
\theta^\star\bigr)\bigr\|_{\infty} >\varepsilon \Bigr) \to0\qquad \mbox{as
} \sigma\to0.
\]

\item[(4)] $\mathbb{P}_{\theta^\star, \sigma}  (\|\sigma^{-1}\nabla
_{S_0} \ell_{Y(\omega)}(\theta^\star)\| <\infty ) =1 $ for
small enough $\sigma$.

\item[(5)] There exists a $p_0\times p_0$ positive definite matrix $\Omega
_{00}$ such that
\[
\forall\varepsilon>0\qquad \mathbb{P}_{\theta^\star, \sigma} \Bigl( \sup_{\theta\in\Theta^\star(\delta)}
\bigl\|\nabla_{S_0, S_0}^2 \ell _{Y(\omega)}(\theta) +
\Omega_{00} \bigr\| >\varepsilon \Bigr) \to0 \qquad \mbox{as } \sigma\to0.
\]
\end{longlist}
\end{ass}

\renewcommand{\theass}{P}
\begin{ass}[(On the prior distribution)]\label{assP}
The $\sigma$-finite measure $\pi(d\theta)$ on $\Theta$ satisfies
the following conditions:
\begin{longlist}[(2)]
\item[(1)] $\int_{\Theta} e^{\ell_y(\theta)/\sigma^2} \pi(d\theta) <
\infty$ for $\mathbb{P}_{\theta^\star,\sigma}$-almost all $y\in
\mathcal{Y}$, for small enough $\sigma$.

\item[(2)] For $\theta\in\Theta^\star(\delta)$, there exists $p(\theta
)\geq0$ such that $\pi(d\theta) = p(\theta) \,d\theta$.

\item[(3)] There exist $C_{\pi} >0$ and $\bolds\alpha_j >0$
for $j\in S_1\cup S_0^\star$, independent of $\sigma$, and there
exists $\Delta_{\pi}= \Delta_{\pi}(\delta) \geq0$, such that
$\Delta_{\pi} \to0$ as $\sigma\to0$ and for $\theta\in\Theta
^\star(\delta)$,
\[
C_{\pi} (1- \Delta_{\pi}) \leq p(\theta) \times\prod
_{j \in S_1
\cup S_0^\star} \theta_j^{-(\bolds\alpha_j-1)} \leq
C_{\pi} (1+ \Delta _{\pi}). %
\]
Denote $\alpha_0 = \bolds\alpha_{S_0^\star}$, $\alpha_1 = \bolds
\alpha_{S_1}$.
\end{longlist}
\end{ass}

\renewcommand{\theass}{L}
\begin{ass}\label{assL}
Assume $\mathbb{P}_{\theta^\star, \sigma}(\Delta_0 ( \delta) \to
0) \to1$ as $\sigma\to0$, where
%
\begin{equation}
\label{eqDefD0} \quad\Delta_0 (\delta) = \sigma^{-p_0 -\sum_{j\in T_0^\star}(\alpha
_{0,j}-1) - 2\sum_{j=1}^{p_1}\alpha_{1,j}} \int
_{\Theta\setminus
\Theta^\star(\delta)} e^{ ( \ell_{Y}(\theta)-\ell_{Y }(\theta
^\star))/\sigma^2} \pi(d\theta).
\end{equation}
\end{ass}

Assumption~\ref{assL} implies consistency of the posterior distribution at a
certain rate,
and it can be written as $\pi( \Theta^\star(\delta)\mid Y) =
1+O_{\mathbb{P}_{\theta^\star, \sigma}}(1)$ as $\sigma\to0$.
Consistency of the posterior is a necessary assumption for the
Bernstein--von Mises theorem [\citet{vdVaart}, Theorem~10.1].
Under Assumption~\ref{assM}, Assumption~\ref{assL} holds if the following condition is satisfied:
%
\begin{eqnarray}\label{eqDefD02}
\sigma^{-p_0 -\sum_{j\in T_0^\star}(\alpha_{0,j}-1) - 2\sum
_{j=1}^{p_1}\alpha_{1,j}} \int_{\Theta\setminus\Theta^\star
(\delta)}e^{- h(\theta)/\sigma^2} \pi(d\theta) \to0
\nonumber\\[-10pt]\\[-10pt]
\eqntext{\mbox {as } \sigma\to0,}
\end{eqnarray}
where the function $h(\theta)$ is such that
\[
\ell^\star(\theta)- \ell^\star\bigl(\theta^\star
\bigr)\leq- h(\theta )\qquad\mbox{for all } \theta\in\Theta\setminus
\Theta^\star (\delta). %
\]

Under Assumption~\ref{assB}, the complement of the polar cone of the set $\Theta
-\theta^\star$ coincides with $\Theta-\theta^\star$ in a small
enough neighbourhood of $0$; this is essential for the analytic
arguments of the paper. This property holds for other polyhedral
boundaries; for affine transformations of the positive orthant this is
trivial, while in general it relies on the fact that $\sigma\to0$.
For a set $\Theta$ that does not satisfy these conditions, the support
of the posterior distribution in the limit may depend on the complement
of the polar cone of $\Theta-\theta^\star$; see also \citet{Shapiro}.

In Assumption~\ref{assS}, we assume uniform convergence in probability of the
derivatives of the scaled log-likelihood at $\theta^\star$ as $\sigma
$ tends to 0, and that the score function of $\theta_{S_0}$ converges
to 0 at rate $\sigma^{-1}$.

In Assumption~\ref{assP}, we assume that the posterior distribution is proper
but we do not assume that the prior measure itself is proper. Neither
do we assume that $p(\theta)$ is finite and bounded away from 0 on the
boundary of the parameter space,
that is, that $\bolds\alpha_j=1$ for all $j$, which is the assumption of
the BvM theorem.
In particular, the log cosh Markov random field prior distribution that
was discussed in Section~\ref{secMotivation} for the motivating
example, satisfies these conditions with $\bolds\alpha_j=1$ for all
$j\in
S_1\cup S_0^\star$. Other improper priors such as the Jeffreys prior
for a Poisson likelihood, as well as the conjugate Gamma prior and Beta
prior conjugate to a binomial likelihood, satisfy this assumption; see
examples in Section~\ref{secBvMexamples}.

\subsection{The main result}\label{secTheoremBvM1}

Before presenting the main result, we state two preliminary lemmas.
Firstly, we show that the elements $\theta^\star_{S_1}$ are on\vspace*{1pt} the
boundary of $\Theta$, and secondly, we study properties of the
derivatives of $\ell^\star(\theta)$.
%
\begin{lemma}\label{lemopt} If Assumption~\ref{assM} in Section~\ref{secLik}
and Assumption~\ref{assB} in Section~\ref{secQuadraticApprox} hold, then
$\theta^\star_{S_1} =0$ and vector $\nabla_{S_1} \ell^\star(\theta
^\star)$ has negative coordinates.

If also for any $\varepsilon>0$, $\mathbb{P}_{\theta^\star, \sigma
} (\|\nabla^2_{S_0,S_0} \ell_Y(\theta^\star) - \nabla
^2_{S_0,S_0} \ell^\star(\theta^\star)\| >\varepsilon ) \to0$
as $\sigma\to0$, then the matrix $\Omega_{00}=-\nabla^2_{S_0, S_0}
\ell^\star(\theta^\star)$ is positive semi-definite.
\end{lemma}
This lemma follows from standard optimality conditions [e.g.,
Proposition~2.1.2 in \citet{Bertsekas}].

Define the following scaling transform ${\mathcal S}= {\mathcal
S}_{\sigma}\dvtx  \Theta-\theta^\star\to\mathbb{R}^{p_0} \times
\mathbb{R}_+^{p_1}$:
%
\begin{equation}
\label{defScaledTransform} {\mathcal S} \bigl(\theta- \theta^\star\bigr) =
D_{\sigma}^{-1} U \bigl(\theta- \theta^\star\bigr),
\end{equation}
where $D_{\sigma} = \operatorname{diag}( \sigma I_{p_0}, \sigma^2
I_{p_1} )$
and $U = (U_0^T\dvtx  U_1^T )^T$ is defined in Section~\ref{secQuadraticApprox}.
The two subsets of coordinates are scaled differently; we are
considering $(\theta_{S_0}-\theta^\star_{S_0})/\sigma$ and $(\theta
_{S_1}-\theta^\star_{S_1})/\sigma^2$.
In the next lemma we study the image of $\Theta^\star(\delta)$
defined by (\ref{eqDefBdelta}) under this transformation, in the limit.

\begin{lemma}\label{lemLimitX} Let Assumption~\ref{assB} in Section~\ref
{secQuadraticApprox} hold, and take $\delta_0$ and $\delta_1$ such
that $\delta_0 \leq c_0$, $ \delta_1\leq c_1$, $\delta_0/\sigma\to
\infty$, and $\delta_1/\sigma^2\to\infty$ as $\sigma\to0$. Then,
\[
\liminf_{\sigma\to0} {\mathcal S}_{\sigma} \bigl(
\Theta^\star(\delta )-\theta^\star\bigr) = \mathbb{R}^{p_0-p^\star_0}
\times\mathbb {R}_+^{p^\star_0+p_1}, %
\]
the $\liminf$ being in the sense of \citet{Shapiro}.
\end{lemma}
The proof of the lemma is given in Appendix~\ref{secProofAux}.

The limit of the posterior distribution is described by the following
parameters: $\alpha_0 = \bolds\alpha_{S_0^\star}$ and $\alpha_1 =
\bolds\alpha_{S_1}$ defined in Assumption~\ref{assP}, $\Omega_{00}$ defined in
Assumption~\ref{assS}, and $a_0(\omega) $ and $a_1 $ defined by
%
\begin{equation}
\label{eqParamDef} a_0(\omega)= \sigma^{-1}
\Omega_{00}^{-1} \nabla_{S_0} \ell _{Y(\omega)}
\bigl(\theta^\star\bigr),\qquad a_1 = -\nabla_{S_1}
\ell^\star\bigl(\theta^\star\bigr).
\end{equation}
The vector $a_1$ has positive coordinates, which follows from
Lemma~\ref{lemopt}.
The matrix $\sigma^{-2}\Omega_{00}$ is an analogue of the Fisher
information for $\theta_{S_0}$.

In the theorem below, which is an analogue of the Bernstein--von Mises
theorem, we claim that under the stated assumptions, the posterior
distribution of ${\mathcal{S}}(\theta-\theta^\star)$, $\mathbb
{P}_{ {\mathcal S}(\theta-\theta^\star)\mid Y}$, converges to a
finite limit.

\begin{theorem}\label{thPostApproxBvM}
Consider the Bayesian model defined in Section~\ref{secSetup} under
Assumption~\ref{assM} and such that Assumptions~\ref{assB},~\ref{assS},~\ref{assP} and~\ref{assL}
hold.

Define a random probability measure on $\mathcal{V}_0 \times\mathbb
{R}^{p_1}_+$, with $v=(v_0, v_1)$:
\[
\mu^\star(\omega) (dv) = \mathcal{PTN}_{p_0}
\bigl(dv_0; a_0(\omega), \Omega_{00}^{-1},
p_0^\star,\alpha_0\bigr) \times
\Gamma_{p_1}(dv_1; \alpha_1,
a_1), %
\]
where $\mathcal{V}_0 = \mathbb{R}^{p_0-p^\star_0}\times\mathbb
{R}^{p^\star_0}_+$,
$\mathcal{PTN}_{p_0} (dv_0; a_0, \Omega_{00}^{-1}, p_0^\star, \alpha
_0)$ is the polynomially-tilted truncated Gaussian distribution defined
by (\ref{eqModifGauss}), and $\Gamma_{p_1} (\cdot; \alpha_1,
a_1 )$ is the probability measure of a $p_1$-dimensional vector
$\xi$ with independent coordinates $\xi_i \sim\Gamma(\alpha
_{1,i},a_{1,i})$.

Then, with transform ${\mathcal S}$ defined by (\ref
{defScaledTransform}), as $\sigma\to0$,
\[
\forall\varepsilon>0\qquad \mathbb{P}_{\theta^\star, \sigma} \bigl( \bigl\| \mathbb{P}_{ {\mathcal S}(\theta-\theta^\star)\mid Y}
- \mu ^{\star} \bigr\|_{\TV} >\varepsilon \bigr) {\to} 0. %
\]
\end{theorem}
The proof is given in Appendix~\ref{secTheProof}. If $\theta^\star$
is an interior point, then $p_1=p^\star_0=0$,
the additional factor in the definition of $\mu^\star$ disappears,
and the limit is Gaussian, as in the classical Bernstein--von Mises theorem.

Assumptions~\ref{assM} and~\ref{assS} imply that the log-likelihood can be approximated
quadratically with respect to the parameter $\theta_{S_0}$ (which
includes $\theta_{S_0^\star}$ where the ``true'' parameter is on the
boundary of the parameter space) but not with respect to~$\theta
_{S_1}$. This is related to the LAN property [\citet{LeCamYang}].
In particular, the rate of convergence for $\theta_{S_0^\star}$ is
still $\sigma^{-1}$, and the limit\vspace*{2pt} of the rescaled posterior is
truncated Gaussian, possibly modified by the behaviour of the local
prior density on the boundary, whereas for $\theta_{S_1}$ the rate of
convergence is faster ($\sigma^{-2}$ instead of $\sigma^{-1}$),
$\theta_{S_1}$ is asymptotically independent of $\theta_{S_0}$ given
data, and its limiting distribution is Gamma. See examples in
Section~\ref{secBvMexamples}.

We shall see in Section~\ref{secBvMexamples} that in a number of
models parameter components on the boundary can only be either all
regular or all nonregular. However, in the motivating SPECT example,
both types of boundary behaviour can occur. Hence the chosen prior,
satisfying Assumption~\ref{assP} with $\bolds\alpha_j=1$ for all $j\in S_1\cup
S_0^\star$, results in asymptotically efficient inference for the
regular parameters.

\begin{remark} The key property of the
posterior distribution, when the true parameter is on the boundary, is
that the gradient of the log-likelihood at this point does not vanish
asymptotically. Thus in some directions the leading term at the Taylor
expansion of log posterior density is linear rather than quadratic, as
would be the case when $\theta^\star$ is an interior point.
If the local prior density at $\theta^\star$ is bounded away from 0
and infinity, then the limit of the posterior in these directions is an
exponential distribution; if the local prior density has an additional
polynomial term in a neighbourhood of $\theta^\star_j =0$, then the
limit is a Gamma distribution.
\end{remark}

If the prior density behaves like a positive constant on the boundary
or the ``regular'' part of the parameter is not on the boundary, then
the limiting distribution $\mu^\star(\omega)$ has a simple form.
%
\begin{corollary}\label{thPostApproxBvMGamma}
Assume that Assumption~\ref{assP} is satisfied with $\alpha_{0,j} =1$ for \mbox{$j\in
T_0^\star$}, or the set $ T_0^\star$ is empty (i.e., $p^\star_0=0$).
Then, under the conditions of Theorem~\ref{thPostApproxBvM}, the
limiting probability measure $\mu^\star(\omega)$ on $\mathcal{V}_0
\times\mathbb{R}_+^{p_1}$ is defined by
\[
\mu^\star(\omega) (dv) = {\mathcal{TN}}_{p_0}
\bigl(dv_0; a_0(\omega ), \Omega_{00}^{-1}
\bigr) \times \bigotimes_{i=1}^{p_1} \Gamma
(dv_{1,i}; \alpha_{1,i}, a_{1,i}), %
\]
where ${\mathcal{TN}}_{p_0} (dv_0; a_0(\omega), \Omega
_{00}^{-1} ) $ denotes the Gaussian distribution truncated to
$\mathcal{V}_0$ and normalised to be a probability measure.

In particular, if the prior distribution behaves as a constant in a
neighbourhood of $\theta^\star$ ($\alpha_{1,j} =1$ for all $j$),
then the limit of $\theta_{S_1}/\sigma^2$ is multivariate exponential.
\end{corollary}

\subsection{Efficiency of inference for ``nonregular'' parameters}

We can see that for $\theta_{S_0\setminus S_0^\star}$ the standard
Bernstein--von Mises theorem holds under the assumption that the prior
density in the neighbourhood of $\theta_{S_0\setminus S_0^\star}$ is
bounded away from 0 and infinity, a standard assumption of the BvM
theorem. Thus inference for $\theta_{S_0\setminus S_0^\star}$ is
asymptotically independent of the prior and is asymptotically
equivalent to efficient frequentist inference.

However, inference for $\theta_{S_1}$ is different. The first key
difference is that there is no need to require a similar assumption on
the prior distribution: even if the local prior density tends to
infinity or zero (both at a polynomial rate) on the boundary, for
i.i.d. observations with $\sigma^2=1/n$, Bayesian inference is still
consistent, at a rate faster than the parametric $\sqrt{n}$ rate. The
second difference is that the limit of the rescaled and recentred
posterior distribution for $\theta^\star_{S_1}$ is not random (i.e.,
does not depend on $\omega$).
These two properties lead to the third important difference which is
the formulation of efficiency of the estimation procedure for these
``nonregular'' parameters. This point is elaborated below.

Consider the case where $p=1$ and $\theta^\star$ is on the boundary
(i.e., $\theta^\star=0$) with $\nabla\ell^\star(\theta^\star
)\neq0$. If the prior density at $\theta^\star$ is not bounded away
from 0 and infinity, the limit of the posterior distribution depends on
the behaviour of the prior distribution on the boundary via exponent
$\alpha$ ($\bolds\alpha_j$ with $j=1$). This exponent is a construct of
the statistician and does not depend on the data or its model and can
be chosen freely. If $\alpha> 1$, then the prior density at the true
value $\theta^\star$ is 0, and if $\alpha<1$, the local prior
density of $\theta$ tends to infinity as $\theta\to\theta^\star$.
The length of the asymptotic posterior credible interval for $\theta$
decreases to 0 as $\alpha\to0$ (see Examples~\ref{exPoisson} and
\ref{exBinom} in Section~\ref{secBvMexamples}); hence it is
possible to recover the true value on the boundary as precisely as
desired, up to the error of approximation of the posterior distribution
by its limit (an upper bound on that is presented in Proposition~\ref
{propBvMnonasympt}). Note that for the Poisson and Binomial
distributions discussed in Examples~\ref{exPoisson} and~\ref{exBinom}, the Jeffreys prior satisfies Assumption~\ref{assP} with $\alpha=1/2$.
This property raises questions about the formulation of efficiency in
this case, as, from a theoretical perspective, there appears to be no
lower bound on the length of the credible interval as in the regular case.

\subsection{Nonasymptotic upper bound}

We also state a nonasymptotic bound on the distance between the
posterior distribution of the rescaled parameter and its limit.

\begin{proposition}\label{propBvMnonasympt}
Assume that the following conditions hold for $\delta_0$ and $\delta
_1$ and for some $\delta_{*0}>0$, $\delta_{*1}>0$ that may depend on
$\delta_0$ and $\delta_1$:
%
\begin{eqnarray}
\label{eqNonAsympCond}
\delta_{*1} &<& a_{\min},\qquad
\delta_{*0}< \lambda_{\min}(\Omega _{00}),\qquad
\delta_0< \bigl\|\theta^\star_{S_0}\bigr\|,
\nonumber\\[-8pt]\\[-8pt] \nonumber
\delta_0 &\leq& c_0,\qquad
\delta_1 \leq c_1,
\end{eqnarray}
where $a_{\min} = \min_{j} a_{1,j}$, $\lambda_{\min}(\Omega_{00})$
is the smallest eigenvalue of $\Omega_{00}$, and $c_0,c_1$ are
constants from Assumption~\ref{assB}.
Let the assumptions of Theorem~\ref{thPostApproxBvM} hold, and define
the following events:
%
\begin{eqnarray}\label{eqDefEventA}
\mathcal{A}_0 &=& \Bigl\{ \omega\dvtx  \sup_{\theta\in\Theta^\star
(\delta)} \bigl\|
\nabla^2_{S_0, S_0} \ell_{Y(\omega)} (\theta) +
\Omega_{00} \bigr\| \leq\delta_{*0}\Bigr\},
\nonumber
\nonumber\\[-8pt]\\[-8pt]\nonumber
\mathcal{A}_1 &=& \Bigl\{ \omega\dvtx  \sup_{\theta\in\Theta^\star
(\delta)} \bigl\|
\nabla_{S_1} \ell_{Y(\omega)} (\theta) + a_1\bigr\|
_{\infty} \leq\delta_{*1} \Bigr\}.
\end{eqnarray}
Then, on $\mathcal{A}= \mathcal{A}_0 \cap\mathcal{A}_1 \cap\{\|
a_0(\omega)\| < \delta_0/\sigma\}$,
%
\begin{eqnarray} \label{eqNonAsymp}
&& \bigl\| \mathbb{P}_{ {\mathcal S}(\theta-\theta^\star)\mid Y} - \mu ^{\star} \bigr\|_{\TV}\nonumber
\\
&&\qquad \leq 2
\max \biggl\{ C_1 \delta_{*1}, p_1 \max
_{j} \Gamma \biggl( \biggl(\frac{a_{1,j}\delta_1}{\sigma^2},\infty \biggr);
\alpha_{1,j}, 1 \biggr) \biggr\}
\nonumber\\[-8pt]\\[-8pt]
&&\quad\qquad{}+ 2 \max \biggl\{ C_0\delta_{*0}, C_{\alpha_0}
\Gamma \biggl( \biggl( \frac{\lambda_{\min}(\Omega
_{00})}{2} \biggl[\frac{\delta_0}{\sigma} -
\bigl\|a_0(\omega)\bigr\| \biggr]^2, \infty \biggr);
\frac{p_{\alpha0}} 2, 1 \biggr) \biggr\}\hspace*{-10pt}
\nonumber
\\
&&\quad\qquad{}+ C_2 \Delta_{\pi} + C_{\Delta}
\Delta_0(\delta),\nonumber
\end{eqnarray}
where\vspace*{1pt} $p_{\alpha0} = p_0 + \sum_{j\in T_0^\star}(\alpha_{0, j}-1)$
and the constants are defined in the proof.
If also $\delta_{*0} \to0$ and $\delta_{*1} \to0$ as $\sigma\to
0$, then the upper bound in (\ref{eqNonAsymp}) tends to 0.
\end{proposition}
The proof is given in Appendix~\ref{secTheProof}. Note that under the
assumptions of Theorem~\ref{thPostApproxBvM}, $\mathbb{P}_{\theta
^\star,\sigma}(\mathcal{A}) \to1$ as $\delta_{*0} \to0$ and
$\delta_{*1} \to0$.
For the upper bound of the total variation to be small in practical
applications, the dimensions $p_k$ should not be too large compared to
the corresponding rate, the smallest eigenvalue of the precision matrix
$\Omega_{00}$ cannot be too small, that is, that $\lambda_{\min
}(\Omega_{00})\delta_0^2/\sigma^2$ should be large, and that
the combination of parameters $(\alpha_{1,j}, a_{1,j})$ should be such
that value $\delta_1/\sigma^2$ is far in the tail of all
corresponding Gamma distributions. If $\alpha_{1,j}=1$ for all $j$,
this requires that the smallest value $a_{\min}$ of the parameter
$a_1$ should not be too small, that is, $a_{\min}\delta_1/\sigma^2$
should be large.

It is interesting to note that, for each $k=0,1$, if $\delta_{*k}
\asymp\delta_k$, which holds in many cases, the value of $\delta_k$
minimising the local upper bound (the first two lines of the upper
bound) coincides with the upper bound of the Ky Fan distance between
the posterior distribution of $\theta_{S_k}$ and its limit, a point
mass at $\theta^\star_{S_k}$. These are $\delta_0 = C_{\Omega_{00}}
\sigma\sqrt{\log(1/\sigma)}$ and $\delta_1 = C_{a_1} \sigma^2
\log(1/\sigma)$ [\citet{BochkinaIP}].

\section{Examples}\label{secBvMexamples}

We now give examples where the asymptotic posterior distribution
differs from Gaussian.
We start with a rule to verify Assumption~\ref{assL} which applies to
exponential family distributions that we consider below.
%
\begin{lemma}\label{lemAssumptionL} Take $\delta_0, \delta_1 >0$
such that $\delta_0, \delta_1 \to0$, and assume that for any $\theta
\in\Theta\setminus\Theta^\star(\delta)$,
\[
\ell_Y(\theta)-\ell_Y\bigl(\theta^\star\bigr)
\leq-C_{\delta0} \sum_{j\in
S_0} \bigl|
\theta_j-\theta^\star_j\bigr| -C_{\delta1}
\sum_{j\in S_1} \bigl|\theta_j-
\theta^\star_j\bigr| %
\]
for some $C_{\delta0}, C_{\delta1} >0$ with probability close to 1
for small enough $\sigma$,
and that there exist $\alpha_j >0$, $j=1,\ldots,p$, and $C_{\pi0}
>0$ such that for all $\theta\in\Theta$,
\begin{eqnarray*}
\frac{\pi(d\theta)}{d\theta} &\leq& C_{\pi0} \prod_{j\in S_0\dvtx
|\theta_j|<\delta_0/\sqrt{p_0}}
\theta_j^{\bolds\alpha_j-1} \prod_{j\in S_1\dvtx  \theta_j<\delta_1}
\theta_j^{\bolds\alpha_j-1}.
\end{eqnarray*}
If $C_{\delta0}\delta_0/\sigma^2 \to\infty$ and $C_{\delta
1}\delta_1/\sigma^2 \to\infty$,
then
$ \Delta_0 (\delta) \to0$ as $\sigma\to0$ with probability 1, that
is, Assumption~\ref{assL} is satisfied.
\end{lemma}
The proof is given in Appendix~\ref{secProofAux}.

%
\begin{example}[(Poisson likelihood)]\label{exPoisson}
Consider $Y_i \sim\operatorname{Poisson}(\theta)$ independently for
$i=1,\ldots, n$, where the true value is $\theta^\star=0$.
In this case, $\sigma^2 = 1/n$ and $\mathbb{P}(Y_i =0)=1$. Consider
an improper prior for $\theta$ with density $p(\theta) = \theta
^{\alpha-1}$ with some $\alpha>0$; the case $\alpha=1/2$ corresponds
to the Jeffreys prior for parameter $\theta$.
In this case, the exact posterior distribution for $\theta$ is $\Gamma
(\alpha, n)$, that is, $n\theta\mid Y \sim\Gamma(\alpha,1)$ which
agrees with Theorem~\ref{thPostApproxBvM}, and the exact $95\%$
credible interval for $\theta$ is $ [0, \gamma_{\alpha}(0.05)/n
 ]$ where $\gamma_{\alpha}(0.05)$ is the 95\% percentile of the
$\Gamma(\alpha,1)$ distribution. For $\alpha= 1/2$, the credible
interval is $[0, 1.92/n]$, and for $\alpha= 0.05$, it is $[0,0.27/n]$.
By decreasing $\alpha$ to 0, we can construct a credible interval of
arbitrarily small length for fixed $n$, even for $n=1$.
\end{example}

%
\begin{example}[(Binomial distribution)]\label{exBinom}
Consider the problem of estimating the unknown probabilities of
Binomial distributions
$Y_i \sim\operatorname{Bin}(n_i, \theta_i)$ independently,
$i=1,\ldots,p$,
for $\theta_i \in[0,1]$, where the true value $\theta^\star_i$ of
some $\theta_i$ is 0.
We assume that all $\theta^\star_i<1$ (if $\theta^\star_i=1$ for
some $i$, consider $n_i - Y_i$ as data and $1- \theta^\star_i$ as the
corresponding parameter). We study the limit of the posterior
distribution for large $n_i$ for all $i=1,\ldots,p$ such that $n_i/n
\to\omega_i \in(0,1)$ where $n=\sum_{i=1}^p n_i$, and $p$ is fixed.
This situation is not covered by the standard BvM theorem. Consider a
conjugate Beta prior $\theta_i \sim B(\alpha, 1)$ independently, with
some fixed $\alpha>0$. In this case, $\sigma^2 = 1/n$
and, as $n\to\infty$,
\begin{eqnarray*}
\ell^\star(\theta) &=& \lim_{n\to\infty}
\ell_{Y}(\theta) = \sum_{i=1}^p
\omega_i\bigl[ \theta^\star_i \log(
\theta_i) - \bigl(1 - \theta^\star_i\bigr)
\log(1- \theta_i)\bigr].
\end{eqnarray*}
If $\theta^\star_i=0$, the corresponding summand in $\ell^\star
(\theta)$ is $-\omega_i\log(1- \theta_i)$ which is defined for
$\theta_i \in[0,1)$, and
then $\nabla_i \ell^\star(\theta) = \omega_i/(1- \theta_i)$.
In this case, $S_0^\star$ is always empty, $\nabla_i \ell^\star
(\theta^\star) = 0$ for $\theta^\star_i\neq0$ and $\nabla_i \ell
^\star(\theta^\star) = - \omega_i$ for $\theta^\star_i= 0$.
Assumption~\ref{assM} was verified in Example~\ref{exBinomM}, and it is easy to check that
Assumptions~\ref{assS},~\ref{assP} and~\ref{assL} are satisfied (e.g., for $p=1$ and $\theta
^\star_i=0$, conditions of Lemma~\ref{lemAssumptionL} hold with
\mbox{$C_{\delta_1}=1$} and $C_{\pi0}=\alpha$). Therefore, $\Omega_{00} =
\operatorname{diag} ( \omega_i/\{\theta^\star_i (1 - \theta
^\star_i)\},
i\in S_0  )$, $a_1 = (\omega_i, i\in S_1)$, and
$a_{0,i} =(Y_i- n_i \theta^\star_i)/(\sqrt{n} \omega_i)$.
Theorem~\ref{thPostApproxBvM} implies the following asymptotic
approximation of the posterior distribution of $ (\theta_{S_0}, \theta
_{S_1} )$:
\[
( \theta_{S_0}, \theta_{S_1} ) \mid Y \sim\mathcal
{N}_{p_0} \biggl( \theta^\star_{S_0} +
\frac{a_0}{\sqrt{n}}, \frac{1}n \Omega_{00}^{-1}
\biggr) \times\Gamma_{p_1} (\alpha, n a_1). %
\]
Similarly to the Poisson likelihood case (Example~\ref{exPoisson}),
for $\alpha$ close to 0, the approximate credible intervals for
$\theta_i$, $i\in S_1$, are small. This is easy to see from the
marginal $100(1-\beta)\%$ credible intervals which are $ [0,
\gamma_{\alpha}(\beta)/(n \omega_i)  ]$.
\end{example}

\begin{example}[(Mixed effects model)]
Consider a model studied by \citet{LRboundary}:
$Y_{ij}\mid\beta_i \sim\mathcal{N}(\mu+\beta_i, \tau^2)$ where
$\beta_i \sim\mathcal{N}(0,\theta)$
independently, for $i=1,\ldots, n$ and $j=1,\ldots, m$.
Here there are $n$ classes with $m$ elements in each, and the parameter
of interest is the contribution of the classes that is characterised by
the parameter $\theta\in\Theta= [0,\infty)$, where the value
$\theta=0$ corresponds to the absence of the random effects $\beta
_i$. We study the asymptotic concentration of the posterior
distribution of $\theta$ when the number of classes $n$ grows while
the number of class elements $m$ remains fixed. We consider a prior
distribution for $\theta$ with density $p(\theta) \propto\theta
^{\alpha-1} e^{-b\theta}$ for $\alpha>0$ and $b\geq0$, which
includes a case of improper prior distributions when $b=0$. Note that
the inverse Gamma prior with density $p(\theta) \propto\theta
^{-\alpha-1} e^{-b/\theta}$ potentially leads to very slow
convergence, since it has a root of infinite order at 0.

We start with the case $\mu$ and $\tau$ known, so without loss of
generality we fix $\mu=0$ and $\tau=1$.
After integrating out $\beta_i$ we have that
$\widebar{Y}_i= m^{-1} \sum_{i=1}^m Y_{ij} \sim\mathcal{N}(0, \theta
^\star+ 1/m)$, independently, where $\theta^\star$ is the true value
of the parameter $\theta$. If $\theta^\star> 0$, then the model is
regular and the posterior distribution of $\theta$ is asymptotically
Gaussian. Now we consider the case $\theta^\star=0$.
Using the marginal likelihood of $ \bar{y}_i$ given $\theta$ and
taking $\sigma^2 = 1/n$, we have
\[
\ell^\star(\theta) = \lim_{n\to\infty}
\ell_{Y}(\theta) = -\frac{1}{2m( \theta+ 1/m)} - \frac{1}{2} \log(
\theta+1/m) %
\]
since $\mathbb{E}\widebar{Y}_i^2 = \theta^\star+ 1/m = 1/m$, and
Assumption~\ref{assM} is satisfied with $\nabla\ell^\star(\theta^\star) =0$.
It is easy to check that Assumptions~\ref{assB},~\ref{assS},~\ref{assP} and~\ref{assL} are satisfied,
and $\nabla^2 \ell^\star(\theta^\star) =-m^2/2$. Thus, by
Theorem~\ref{thPostApproxBvM}, the approximate posterior distribution
of $\sqrt{n} \theta$ has density
\[
p_{\theta\sqrt{n}}(x\mid y) \approx C_{\alpha, m, a_0} x^{\alpha-1}
e^{-(x- a_0)^2/m^2},\qquad x\geq0 %
\]
with $a_0 =(m/(2 \sqrt{n})) ( n^{-1} \sum_{i=1}^n m\widebar{Y}_i^2 -1 )$.
It is easy to show that the Cramer--Rao lower bound on the variance of
estimators of $\theta$ applies here, even in the case $\theta^\star
=0$. Thus, using a prior with $\alpha< 1$ (i.e., introducing a bias
towards 0) would lead to superefficiency, that is, loss of efficiency
for $\theta^\star\neq0$. In the case $\alpha=1$ the posterior
distribution is Gaussian with the same mean and variance as in the BvM
theorem but truncated to $\theta\geq0$. The length of the credible
interval for $\theta$ in this case is smaller than in the case where
$\theta^\star$ is an interior point.

Now consider the case where parameters $(\mu, \tau^2, \theta)$ are
estimated jointly with a continuous prior for $(\mu, \tau^2)$ whose
density is finite and positive at the true value $(\mu^\star, \tau
^\star)$. Then
\[
-\ell^\star\bigl( \mu,\tau^2, \theta\bigr) =
\frac{ (m-1){\tau^\star
}^2}{2\tau^2} + \frac{ (\mu-\mu^\star)^2 + {\tau^\star}^2/m
}{2\tau^2 ( \theta+ 1/m)} + \frac{\log(\theta+1/m)}{2} + \frac{m
\log(\tau^2)}{2},
\]
since $\mathbb{E}(\widebar{Y}_i-\mu)^2 = {\tau^\star}^2/ m +(\mu-\mu
^\star)^2$ and $\mathbb{E}\sum_{j=1}^m (Y_{ij}- \widebar{Y}_i)^2 =
(m-1){\tau^\star}^2$. The function $ \ell^\star( \mu,\tau^2,
\theta)$ is maximised at $\mu=\mu^\star$, $\tau= \tau^\star$ and
$\theta= \theta^\star=0$, with zero gradient and the negative matrix
of the second order derivatives $\Omega_{00}$ and its inverse (the
covariance matrix) being
\begin{eqnarray*}
\Omega_{00} &=& \pmatrix{ \displaystyle\frac{ m}{{\tau^\star}^2} &0 & 0
\cr
0
& \displaystyle\frac{ m }{2{\tau^\star}^4} & \displaystyle\frac
{m}{2 {\tau^\star}^2}
\cr
0 &
\displaystyle\frac{m}{2 {\tau^\star}^2}& \displaystyle\frac
{m^2} 2},
\\
\Omega_{00}^{-1} &=& \pmatrix{ \displaystyle\frac{{\tau^\star}^2}{ m}
& 0 & 0
\cr
0 & \displaystyle\frac{2{\tau^\star}^4}{ m -1} & \displaystyle -
\frac{2 {\tau^\star}^2}{m(m-1)}
\cr
0 & \displaystyle-\frac{2 {\tau^\star}^2}{m(m-1)} &\displaystyle
\frac{2} {m(m-1)}}.
\end{eqnarray*}

If $\alpha= 1$, then the approximate joint posterior distribution of
$\sqrt{n}(\theta-\theta^\star, \mu-\mu^\star, \tau^2 - {\tau
^\star}^2)$ is Gaussian truncated to $\theta-\theta^\star= \theta
\geq0$ with bias as given\vspace*{1pt} in Theorem~\ref{thPostApproxBvM} and the
covariance matrix $\Omega_{00}^{-1}$ given above. Note that $\theta$
and $\tau^2$ are asymptotically correlated, with correlation $-m^{-1/2}$.
\end{example}
%
\section{Asymptotic behaviour of the posterior distribution for SPECT}\label{secSPECTexample}

\subsection{Approximation of the posterior distribution}\label{secSPECTexamplePost}

Consider the SPECT model defined in Section~\ref{secMotivation}, in
which $\theta^\star$ has some zero coordinates. The assumptions of
Theorem~\ref{thPostApproxBvM} were verified in Examples~\ref{exPoissonM} and \ref{exPoisson}
(Assumptions~\ref{assM},~\ref{assB},~\ref{assS}), and the log cosh Markov random field prior
distribution satisfies Assumption~\ref{assP} with $\alpha_j=1$ for all $j$.
Assumption~\ref{assL} also holds, since the conditions of Lemma~\ref{lemAssumptionL} are satisfied for independent Poisson random
variables with
$
C_{\delta0} = 0.5 \delta_0 (\delta_0+\sqrt{p_0} y^\star_{\min
})^{-1} y^\star_{\min}$, $C_{\delta1} = \min_j(a_{1,j})$,
where $y^\star_{\min} = \min_{j\dvtx  y^\star_j>0}y^\star_j$ with
$y^\star= A \theta^\star$ for small enough $\sigma=1/\sqrt
{\mathcal{T}}$, due to the inequality $\log(1+x) -x\leq- xb/(b+1)$
for $x>b>0$.

For this model, $\nabla\ell^\star(\theta^\star) =-\sum_{i\dvtx
y^\star_i=0} A_i^T$,
which is nonzero if $Z= \{i \in\{1,\ldots,n\}\dvtx  y^\star_i=0\}$ is
not empty. Hence, nonregularity arises from the elements where there
are no detected photons ($y^\star_i =0$) and the likelihood
degenerates: $\mathbb{P}_{\theta^\star} (Y_i=0)=1$ for $i\in Z$ but,
since $A_i \neq0$, this gives us information about those~$\theta_j$
where $A_{ij} \neq0$, that is, on ${S_1}=\{j\dvtx  \theta^\star_j=0$ and
$\sum_{i \in Z} A_{ij}\neq0\}$.

The limiting distribution of $\theta_{S_1}/\sigma^2$ is exponential
with parameter $a_1 = \sum_{i \in Z} A_{i {S_1}}^T$.
The parameter $(\theta_{S_0} - \theta^\star_{S_0})/\sigma$ has
approximately a truncated Gaussian distribution with parameters
\[
\Omega_{00} = A_{\widebar{Z}, S_0}^T \operatorname{diag}
\bigl(1/\bigl[y^\star \bigr]_{\widebar{Z}}\bigr) A_{\widebar{Z}, S_0},
\qquad a_0 = \Omega_{00}^{-1}
A_{\widebar{Z},S_0}^T \widetilde{Y} /\sigma, %
\]
where $\widetilde{Y}$ is a vector with coordinates $Y_i/y^\star_i -1 $
for $i\in\widebar{Z}$. Truncation takes place for parameters $\theta
_{S_0^\star}$ with $S_0^\star=\{j\dvtx  \theta^\star_j=0$ and $\sum_{i \in Z} A_{ij}= 0\}$.

If the vector of Poisson means $y^\star= A\theta^\star$ has only
positive coordinates ($Z$ is empty), the model is regular, and the
posterior distribution of $(\theta-\theta^\star)/\sigma$ is
approximately truncated Gaussian.

\subsection{Practical implications of the approximate posterior}

We will briefly discuss some practical implications of Theorem~\ref
{thPostApproxBvM}.
Well-developed methods for SPECT reconstruction using our model, using
Markov chain Monte Carlo computation, deliver both approximate,
simulation-consistent, posterior means and variances; see \citet
{Weir} for a fully Bayesian reconstruction. The theorem provides
valuable knowledge which can enrich the interpretation of such results,
enabling approximate probabilistic inference.

Inferential questions of real interest, including (a)
quantitative inference about amounts of radio-labelled tracer within
specified regions of interest, or (b)
tests for significance of apparent hot- or cold-spots,\vspace*{1pt}
can be answered using approximate posteriors for \emph{linear
combinations} $w^T\theta$ of parameters, and are particularly amenable
to treatment. Specifically, suppose that for any nonempty set of pixels
\mbox{$R\subseteq\{1,2,\ldots,p\}$}, $w^R$ denotes the vector with elements
$w^R_j=1/|R|$ for $j\in R$, 0~otherwise. Then to deal with case (a) we
can take $w=w^R$ to deliver $w^T\theta$ as the average concentration
of tracer in region $R$, and for case (b) take $w=w^{R_1}-w^{R_2}$ for
the difference in average concentration between regions $R_1$ and~$R_2$.

To construct an approximation of the posterior distribution, we require
estimates of unknown parameters. We use the marginal posterior modes
estimate $\hat{\theta}$, $\hat{\theta}_i=\operatorname
{argmax} p(\theta
_i|y)$, instead of $\theta^\star$, $\hat{y} = A \hat{\theta}$ instead of $y^\star$,
\[
\widehat{S}_1 = \bigl\{j\dvtx  \nabla_j
\ell_y(\hat{\theta}) < 0\bigr\}, \qquad\widehat{Z} = \{i\dvtx
\hat{y}_i =0\}.
\]
A more robust way to estimate ${S_1}$ would be to use $\widehat
{S}_{1,\varepsilon} = \{j\dvtx  \nabla_j \ell_y(\hat{\theta}) <
-\varepsilon\}$ for some small enough $\varepsilon>0$; however, sensitivity
to the choice of $\varepsilon$ would need to be investigated.
Then, the approximate posterior of $z = (\theta- \hat{\theta})$ is
\[
\phi(z) = \prod_{j \in\widehat{S}_1} \bigl[\hat{a}_j/
\bigl(2\sigma^2\bigr)\bigr] \bigl(2\pi\sigma^2
\bigr)^{-p_0/2} \bigl[\det(\widehat{\Omega})\bigr]^{1/2} \exp\bigl
\{ - z_{\widehat{S}_0}^T \widehat{\Omega} z_{\widehat{S}_0}/\bigl(2
\sigma^2\bigr) - z_{\widehat{S}}^T \hat{a}/
\sigma^2\bigr\}, %
\]
where $\widehat{\Omega} = \sum_{i \notin\widehat{Z}}y_i/[\hat
{y}_i]^2 A_{i, \widehat{S}_0} A_{i, \widehat{S}_0}^T$ and $\hat
{a} = \sum_{i\in\widehat{Z}} A_{i, \widehat{S}_1}^T$.

\subsection{Finite sample performance}\label{secSPECThist}

We briefly discuss the extent to which the approximation in Theorem
\ref{thPostApproxBvM} holds true for data on the scale of a real
SPECT study. A formal assessment of this would entail a major study
beyond the scope of this paper, so we present selected results from
analysis of two
data sets based on a SPECT scan of the pelvis of a human subject.

In the first experiment, the matrix $A$ was constructed according to
the model in \citet{Green} and \citet{Weir}, capturing
geometry, attenuation, and radioactive decay for a setup consisting of
64 projections from a 2-dimensional slice through the patient, each
projection yielding an array of 52 photon counts, on a spatial
resolution of 0.57~cm.
The data set was obtained from Bristol Royal Infirmary; the total
photon count was 45,652; individual counts ranged from 0 to 85,
averaging~13.7.
Reconstruction was performed on a $48\times48$ square grid of 0.64~cm
pixels, using the log cosh prior with hyperparameters fixed at $\gamma
=25$ and $\zeta=8$, obtained using a simple MCMC sampler. We employed
20,000 sweeps of a deterministic-raster-scan single-pixel random walk
Metropolis sampler on a square-root scale for $\theta$, chosen to
avoid extremes in acceptance rate at high- and low-spots in the image.

%
\begin{figure}[b]

\includegraphics{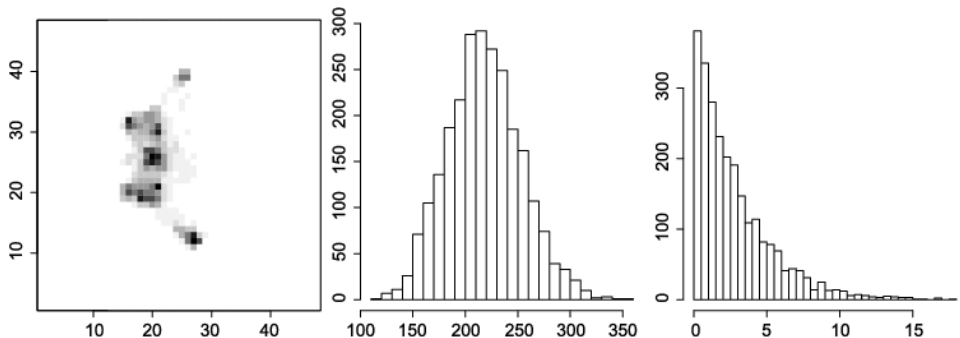}

\caption{Analysis of real SPECT data: posterior mean reconstruction as
a grey-scale image, histogram of marginal posterior for a high-spot
pixel (row 12, column 28),
and the same for a low-spot pixel (row 12, column 31).}
\label{figbrimh1}
\end{figure}

Figure~\ref{figbrimh1} shows selected aspects of this analysis; see
caption for details. Our tentative conclusion is that the marginal
posterior distributions for individual pixels~$\theta_j$ do appear to
be approximately Gaussian in high-spots and approximately exponential
in low-spots, consistent with the theoretical limits presented in
Theorem~\ref{thPostApproxBvM}.

A second experiment was focussed on a more precise and quantitative
assessment of the approximation to the posterior derived in the
previous section. The setup is the same as in the first experiment,
except at half the resolution, so that reconstruction was on a $24
\times24$ grid of 1.28~cm pixels. The corresponding $A$ matrix is now
better-conditioned, and $p$ is only 576, so that manipulation of the
matrices is entirely tractable.
Synthetic data was generated using this $A$ and a ``ground truth''
obtained from an approximate MAP reconstruction from the same real data
set as used above,
yielding photon counts between 0 and 243, totalling 138,310. 50,000
sweeps of the MCMC sampler were used, with prior settings $\gamma
=200$, $\zeta=8$.

\begin{figure}[t]

\includegraphics{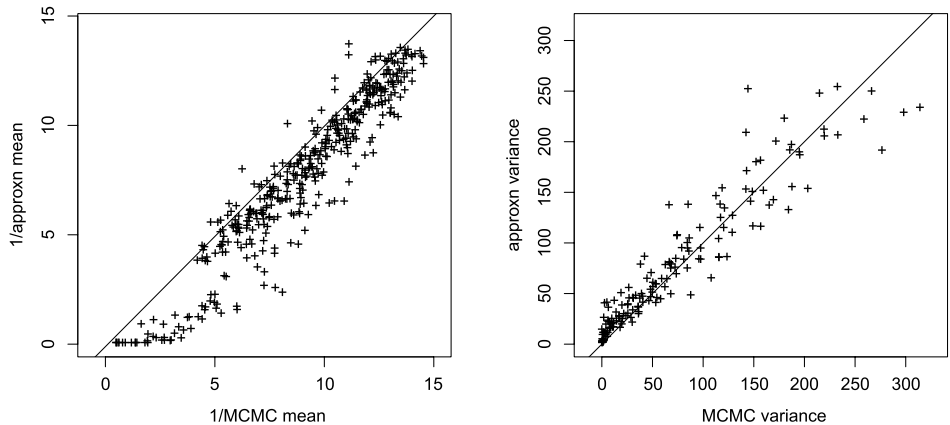}

\caption{Agreement between (left panel) the elements of $\hat{a}$
and the reciprocals of the MCMC-computed posterior means of $\theta$,
for pixels in $\widehat{S}_1$, and also that between (right panel) the
diagonal elements of $\widehat{\Omega}{}^{-1}$ and the posterior
variances of $\theta$ for pixels in $\widehat{S}_0$.}
\label{figapproxnpars}
\end{figure}

Figure~\ref{figapproxnpars} displays the agreement between the
elements of $\hat{a}$ and the reciprocals of the MCMC-computed
posterior means of $\theta$, for pixels in $\widehat{S}_1$, and also
that between the diagonal elements of $\widehat{\Omega}{}^{-1}$ and the
posterior variances of $\theta$ for pixels in~$\widehat{S}_0$.

\begin{figure}

\includegraphics{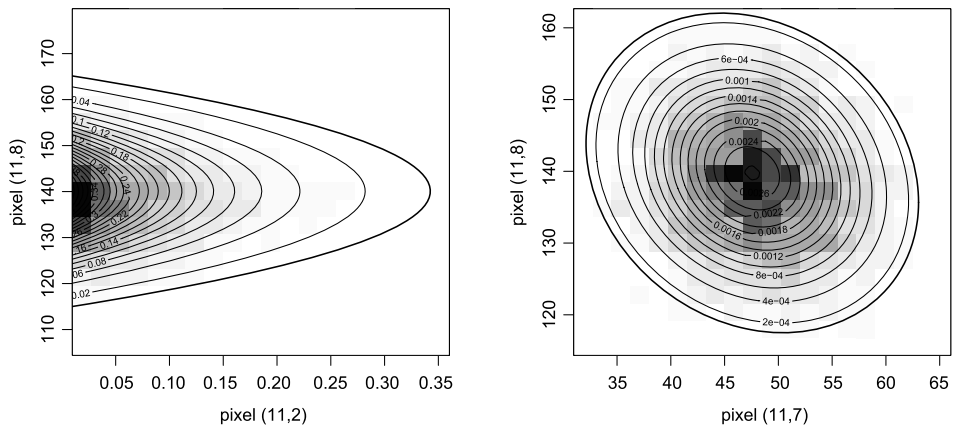}

\caption{Two bivariate\vspace*{1pt} marginals of the posterior, as computed by MCMC
(grey-scale image), and the corresponding approximations (contours). In
the left panel, one pixel is in $\widehat{S}_1$ and one in $\widehat
{S}_0$, so the approximation is Gaussian/exponential; in the right
panel both pixels are from $\widehat{S}_0$, so we have a bivariate
Gaussian. The outermost contour represents the 95\% HPD credible region
based on the approximation.}
\label{figpostcomp}
\end{figure}

Figure~\ref{figpostcomp} displays two bivariate posterior marginals,
computed by MCMC, and the corresponding approximations. In the left
panel, one component is in $\widehat{S}_1$ and one in $\widehat
{S}_0$, so the approximation is Gaussian/exponential; on the right both
components are from $\widehat{S}_0$, so we have a bivariate Gaussian.

We conclude that for this realistic/modest-scale SPECT reconstruction
problem, the small-variance asymptotics\vspace*{1pt} of this paper provide a good
approximation to the posterior, even for $\sigma^2=1$.

\section{Discussion}
When the posterior distribution concentrates on the boundary, we have
shown that the classic Bernstein--von Mises theorem does not hold for
all components.
There are two different types of non-Gaussian component: one, with the
same parametric rate of convergence, is a truncated Gaussian or a
polynomially tilted modification of this if the prior density is not
bounded away from zero and infinity on the boundary, and the second is
a Gamma, with a faster rate of convergence. An interesting property of
the components of the second type is that they are not subject to a
lower bound on efficiency, unlike the ``regular'' and the first-type
boundary components. Under some models with this property, at least
part of the data is observed exactly, so perhaps it should not be an
unexpected phenomenon; see examples of Poisson and Binomial likelihoods
in Section~\ref{secBvMexamples}. This property is quite remarkable:
in principle, it allows the recovery of the unknown parameter on the
boundary with an arbitrarily small precision (particularly in the case
there is no approximation error), by choosing an appropriate prior
distribution, without losing asymptotic efficiency if the parameter is
not on the boundary. This property is related to convergence in
finitely-many steps of the projected gradient method for a sharp
minimum for a noise-free function [\citet{Polyak},
Theorem~1, page~182; thanks to Alexandre Tsybakov for bringing this to
our attention].

A related but different problem involves a nonregular model where the
density of the observations has one or more jumps at a point that
depends on the unknown parameter, for example, $Y_i\sim U[0,\theta]$,
$i=1,\ldots,n$, independently. This type of problem has been
extensively studied from both frequentist and Bayesian perspectives
[\citet{IH81,GGS94,GS95,GGS95,ChernozhukovHong,HiranoPorter}].
In the problem treated in this paper, the rate of convergence of the
posterior distribution of the unknown nonregular parameter as a
function of $n$ is the same as in this case where the unknown parameter
controls the positions of jumps, faster than the standard parametric
rate. However, there is a crucial difference: in the former case, the
posterior distribution has a data-dependent random shift, whereas in
the latter case there is no such shift.

The nonasymptotic version of the main result shows that other
parameters of the model can also affect convergence in practice, such
as the smallest eigenvalues of the precision matrices in the $\mathcal
{PTN}$ part of the limit and the smallest parameter of the scale of the
Gamma distributions.

It is easy to verify that Theorem~\ref{thPostApproxBvM} derived here applies also to
misspecified models, with $\mathbb{P}_{\theta^\star,\sigma}$ being
replaced by the true distribution of $\mathbf{Y}$ and $\theta^\star$
defined as the unique maximum of $\ell^\star(\theta)$ as in
Assumption~\ref{assM}. This will be discussed elsewhere.

An interesting direction for future work is to study both the behaviour
of the posterior distribution, and the question of optimal prior
specification, in a framework where the spatial resolution is
infinitely refined, placing smoothness class constraints on $\theta
^\star$.

\begin{appendix}\label{append}
\section*{Appendix: Proofs}\label{secproofs}
\subsection{Proof of the main result}\label{secTheProof}
We start with a lemma.

%
\begin{lemma}\label{lemApproxNonrandom} Consider the function $\ell
_Y(\theta)$ defined in Section~\ref{secLik} and assume that
Assumptions~\ref{assM},~\ref{assB} and~\ref{assS} hold.
Then, on the event $\mathcal{A}_0 \cap\mathcal{A}_1$ defined by
(\ref{eqDefEventA}) with some $\delta_{*0}, \delta_{*1} >0$, for
$\theta\in\Theta^\star(\delta)$,
\begin{eqnarray*}
&& \ell_Y(\theta) - \ell_Y\bigl(\theta^\star
\bigr)
\\
&&\qquad \geq\bigl(\theta_{S_0}-\theta ^\star_{S_0}
\bigr)^T \nabla_{S_0} \ell_Y\bigl(
\theta^\star\bigr) - \bigl(\theta _{S_0}-\theta^\star_{S_0}
\bigr)^T \widetilde\Omega_{00} \bigl(\theta
_{S_0}-\theta^\star_{S_0}\bigr)/2 -
\tilde{a}^T \theta_{S_1},
\\
&& \ell_Y(\theta) - \ell_Y\bigl(\theta^\star
\bigr)
\\
&&\qquad \leq\bigl(\theta_{S_0}-\theta ^\star_{S_0}
\bigr)^T \nabla_{S_0} \ell_Y\bigl(
\theta^\star\bigr) - \bigl(\theta _{S_0}-\theta^\star_{S_0}
\bigr)^T \widebar\Omega_{00} \bigl(\theta_{S_0}-
\theta ^\star_{S_0}\bigr)/2 -\bar{a}^T
\theta_{S_1},
\end{eqnarray*}
where
\begin{eqnarray*}
\widetilde{\Omega}_{00} &=& \Omega_{00} +
\delta_{*0} I_{p_0},\qquad \widebar{\Omega}_{00} =
\Omega_{00} - \delta_{*0} I_{p_0},
\\
\tilde{a} &=& a_1 + \delta_{*1} \mathbf{1}_{p_1},
\qquad\bar{a} = a_1- \delta_{*1} \mathbf{1}_{p_1}.
\end{eqnarray*}
Here $\mathbf{1}_{p_1} = (1,\ldots, 1)^T$---a vector of length $p_1$,
and $I_{p_0}$ is $p_0 \times p_0$ identity matrix.
\end{lemma}

\begin{pf}
Applying the Taylor expansion of $\ell_Y(\theta)$ as a function of
$\theta_{S_1}$ at point~$\theta^\star_{S_1}$,
and then\vspace*{1pt} expanding $\ell_Y(\tilde\theta)$ where $\tilde\theta
_{S_0} = \theta_{S_0}$ and $\tilde\theta_{S_1} = \theta^\star
_{S_1}$, as a function of $\theta_{S_0}$ at point $\theta^\star
_{S_0}$, for some $\theta_{c0}, \theta_{c1} \in\Theta^\star(\delta
)$, we have
\begin{eqnarray*}
\ell_Y(\theta)- \ell_Y\bigl(\theta^\star
\bigr) &=& \bigl(\theta_{S_1} -\theta ^\star_{S_1}
\bigr)^T \nabla_{S_1} \ell_Y(
\theta_{c1}) + \bigl(\theta_{S_0} -\theta^\star_{S_0}
\bigr)^T \nabla_{S_0} \ell_Y(\theta)
\\
&&{} + \bigl(\theta_{S_0} -\theta^\star_{S_0}
\bigr)^T \nabla_{S_0, S_0} \ell _Y(
\theta_{c0}) \bigl(\theta_{S_0} -\theta^\star_{S_0}
\bigr)/2.
\end{eqnarray*}
Applying the bounds defining events $\mathcal{A}_0$ and $\mathcal
{A}_1$ to $\nabla_{S_1} \ell_Y(\theta_{c1})$ and\break  $\nabla_{S_0, S_0}
\ell_Y(\theta_{c0})$,
and using that $\theta _{S_1}-\theta^\star_{S_1}=\theta_{S_1}$
is a vector with nonnegative components,
we have
\begin{eqnarray*}
\ell_Y(\theta) - \ell_Y\bigl(\theta^\star
\bigr) &\leq& \bigl(\theta_{S_1} -\theta^\star_{S_1}
\bigr)^T [ - a_1 + \delta_{*1}
1_{|{S_1}|}] +\bigl(\theta _{S_0} -\theta^\star_{S_0}
\bigr)^T \nabla_{S_0} \ell_Y(\theta)
\\
&&{} + \bigl(\theta_{S_0} -\theta^\star_{S_0}
\bigr)^T [-\Omega_{00}+ \delta _{*0}
I_{|S_0|}] \bigl(\theta_{S_0} -\theta^\star_{S_0}
\bigr)/2,
\end{eqnarray*}
and hence the first statement of the lemma.
Applying the inequalities on the events~$\mathcal{A}_k$ as lower
bounds, we obtain the second statement of the lemma.
\end{pf}



\begin{pf*}{Proof of Theorem~\ref{thPostApproxBvM}}
Denote $v =(v_0^T, v_1^T)^T = D^{-1} U(\theta-\theta^\star)$ where
$v_0 = (\theta_{S_0} - \theta^\star_{S_0})/\sigma$ and $v_1 =
(\theta_{{S_1}} - \theta^\star_{{S_1}})/\sigma^2$; the Jacobian of
this change of variables is $\sigma^{p_0+2p_1}$.
The image of $\Theta^\star(\delta)$
under this transform is
\[
B_R= B_{2}(0, R_0) \times[0,
R_1)^{p_1} \cap D_{\sigma}^{-1} U \bigl(
\Theta-\theta^\star\bigr),
\]
with $R_0 = \delta_0/\sigma$ and $R_1 = \delta_1 /\sigma^2$. Under
Assumptions~\ref{assB} and~\ref{assS}, the conditions of Lemma~\ref{lemLimitX} hold,
which implies that if $\| \theta^\star_{S_0}\|\geq\delta_0$ and
$\delta_k \leq c_k$, $B_R= [B_{2, p_0}(0, R_0)\cap\mathcal{V}_0 ]
\times[0, R_1]^{p_1}$ where $\mathcal{V}_0 =\mathbb{R}^{p_0-p^\star
_0} \times\mathbb{R}_+^{p^\star_0}$, and the set $B_R$ becomes
$\mathcal{V}^\star=\mathcal{V}_0 \times\mathbb{R}_+^{ p_1}$ as
$\sigma\to0$.

The triangle inequality for the total variation norm gives
%
\begin{eqnarray}\label{eqTriangleTV}
&& \bigl\| \mathbb{P}_{{\mathcal S}(\theta- \theta^\star) \mid Y} - \mu ^\star\bigr\|_{\TV}\nonumber
\\
&&\qquad \leq
\bigl\|\mathbb{P}_{{\mathcal S}(\theta- \theta^\star) \mid Y} \mathbf {1}_{B_R}- \mu^\star
\mathbf{1}_{B_R}\bigr\|_{\TV}
\\
&&\quad\qquad{}+ \bigl\| \mu^\star\mathbf{1}_{B_R} - \mu^\star
\bigr\|_{\TV} + \bigl\|\mathbb {P}_{{\mathcal S}(\theta- \theta^\star) \mid Y}\mathbf{1}_{B_R} -
\mathbb{P}_{{\mathcal S}(\theta- \theta^\star) \mid Y} \bigr\|_{\TV}, \nonumber
\end{eqnarray}
where the balls $B_R$ are defined above. Here $\mu\mathbf{1}_{B_R}$
is a probability measure $\mu$ truncated to $B_R$ and normalised to be
a probability measure.
If the measure $\mu_1$ is absolutely continuous with respect to
measure $\mu_2$, with density $f$, the total variation norm can be
written as
\[
\|\mu_1 - \mu_2\|_{\TV} = 2\int
_{\Theta} (f - 1)_+ \,d\mu_2, %
\]
where $(x)_+ = \max(x,0)$ [\citet{vdVaart}]. This can be used in
each of the summands in the upper bound (\ref{eqTriangleTV}).

In this proof we will use $\alpha= (\alpha_0, \alpha_1)$, for
simplicity of notation.

Define the measure $\mu(dv; a_1, \alpha, b, \Sigma)$ for $v =
(v_0^T, v_1^T)^T$, $v_0\in\mathbb{R}^{p_0-p_0^\star} \times
[0,\infty)^{p_0^\star}$ and $v_1 \in[0,\infty)^{p_1}$, by
%
\begin{equation}
\label{eqMeasureMu} \frac{\mu(dv; a_1, \alpha, b, \Sigma)}{dv} = \prod_{j\in T_0^{\star}\cup T_1}
v_{j}^{\alpha_j-1} e^{ - a_1^T
v_1 - v_0^T \Sigma v_0 /2 + v_0^T b},
\end{equation}
where $T_0^{\star} =\{p_0-p^\star_0+1,\ldots,p_0\}$, $T_1 = \{
p_0+1,\ldots,p\}$, $a_1 \in(0,\infty)^{p_1}$, $b \in\mathbb
{R}^{p_0}$, $\alpha= (\alpha_j)_{j \in T_0^{\star} \cup T_1}\in
(0,\infty)^{p^\star_0+p_1}$, and $\Sigma$ is a $p_0\times p_0$
positive definite matrix.

We start with the first term in (\ref{eqTriangleTV}).
By Lemma~\ref{lemApproxNonrandom}, on the event $\mathcal{A}_0\cap
\mathcal{A}_1$ defined by~(\ref{eqDefEventA}), for any measurable
$\mathcal{B}\subseteq\Theta^\star(\delta)$, with $\mathcal{B}_v =
D_{\sigma}^{-1} U (\mathcal{B}-\theta^\star) \subseteq B_R$, we have
\begin{eqnarray*}
&&\int_{ \mathcal{B}} \exp \bigl\{\bigl[\ell_Y(\theta) -
\ell_Y\bigl(\theta ^\star\bigr)\bigr]/\sigma^2
\bigr\} \pi(d\theta)
\\
&&\qquad \geq J_\sigma C_{\pi} (1-\Delta_{\pi})
\\
&&\quad\qquad{}\times
\int_{ \mathcal{B}_v} \prod_{j\in T_0^{\star}\cup T_1}
v_{j}^{\alpha_j-1} \exp \bigl\{v_0^T
\nabla_{S_0} \ell_Y\bigl(\theta ^\star\bigr)/
\sigma- \bigl\|\widetilde{\Omega}_{00}^{1/2}v_{0}
\bigr\|^2/2 - \tilde{a}^T v_1 \bigr\} \,dv
\\
&&\qquad = J_{\sigma} C_{\pi} (1 - \Delta_{\pi}) \mu\bigl(
\mathcal{B}_v; \tilde{a},\alpha, \nabla_{S_0}
\ell_Y\bigl(\theta^\star\bigr)/\sigma, \widetilde{
\Omega}_{00}\bigr),
\end{eqnarray*}
where $J_\sigma= \sigma^{p_0 -p^\star_0+\sum_{j\in T_0^{\star}}
\alpha_{0,j} + 2\sum_{j=1}^{p_1}\alpha_{1,j}}$, and the measure $\mu
(dv; a_1, \alpha, b, \Sigma)$ is defined by (\ref{eqMeasureMu}).
Similarly, using Lemma~\ref{lemApproxNonrandom}, we obtain an upper
bound on the event~$\mathcal{A}_0 \cap\mathcal{A}_1$,
\begin{eqnarray*}
&&\int_{ \mathcal{B}} \exp \bigl\{ \bigl[\ell_Y(\theta)
- \ell _Y\bigl(\theta^\star\bigr)\bigr]/
\sigma^2 \bigr\} \pi(d\theta)
\\
&&\qquad \leq J_\sigma C_{\pi}
(1+\Delta_{\pi})
\\
&&\quad\qquad{} \times\int_{\mathcal{B}_v} \prod_{j \in T_0^\star\cup T_1}
v_{j}^{\alpha_{j}-1} \exp \bigl\{v_0^T
\nabla_{S_0} \ell_Y\bigl(\theta ^\star\bigr)/
\sigma- \bigl\|\widebar{\Omega}_{00}^{1/2}v_{0}
\bigr\|^2/2 - \bar{a}^T v_1 \bigr\} \,dv
\\
&&\qquad = J_\sigma C_{\pi} (1+\Delta_{\pi}) \mu\bigl(
\mathcal{B}_v; \bar{a}, \alpha, \nabla_{S_0}
\ell_Y\bigl(\theta^\star\bigr)/\sigma,\widebar{
\Omega}_{00}\bigr).
\end{eqnarray*}
To simplify the notation, denote $a_0 = \Omega_{00}^{-1} \nabla_{S_0}
\ell_Y(\theta^\star)/\sigma$ and
\begin{eqnarray*}
\bar{\mu} (dv) &=& \mu(dv; \bar{a}, \alpha, \Omega_{00}a_0,
\widebar {\Omega}_{00}),
\\
\tilde{\mu} (dv) &=& \mu(dv;
\tilde{a},\alpha, \Omega _{00}a_0, \widetilde{
\Omega}_{00}).
\end{eqnarray*}
The measure $ \tilde{\mu}$ is finite since $a_0 = \nabla_{S_0}
\ell_Y(\theta^\star)/\sigma$ is finite with high probability due to
Assumption~\ref{assS}(4), and all its other parameters are positive or positive definite.
The measure $\bar\mu$ is finite if $\delta_{*1} < \min_j a_{1,j}$
and $\delta_{*0} < \lambda_{\min}(\Omega_{00} )$. These conditions
hold if $\delta_{*0}, \delta_{*1}$ are small enough which is possible
due to Assumption~\ref{assS}.

For $\mathcal{B}_v =\mathcal{B}_1\times B_{\infty}(0, r_1)$ for some
$\mathcal{B}_1\subset\mathcal{V}_0$ and $r_1 \in(0, R_1)$, we have
\begin{eqnarray*}
\mu\bigl( \mathcal{V}^\star; a_1, \alpha, b, \Sigma\bigr)
&=& \prod_{i=1}^{p_1} \bigl[a_{1,i}^{-\alpha_{1,i}}
\Gamma(\alpha_{1,i})\bigr] \int_{\mathcal{V}_0} \prod
_{j\in T_0^\star} v_{0, j}^{\alpha_{0,j}-1}
e^{ - v_0^T \Sigma v_0 /2 + v_0^T b} \,d v_0,
\\
\frac{\mu( \mathcal{B}_v; a_1, \alpha, b, \Sigma) }{ \mu(
\mathcal{V}^\star; a_1, \alpha, b, \Sigma)} &=& \mathcal {PTN}_{p_0}
\bigl(\mathcal{B}_1; \Sigma^{-1} b, \Sigma^{-1},
p_0^\star,\alpha_0\bigr) \prod
_{j=1}^{p_1} \Gamma\bigl((0,r_1);
\alpha_{1,j}, a_{1,j}\bigr),
\end{eqnarray*}
where the probability measure $\mathcal{PTN}_{p_0}(\cdot; b, \Omega
_{00}^{-1}, p_0^\star,\alpha_0)$ is defined by (\ref
{eqModifGauss}), and~$\Gamma( \cdot; \alpha_{1,j}, a_{1,j} )$ is
the probability measure associated with distribution $\Gamma(\alpha
_{1,j}, a_{1,j})$.

Hence, the posterior density of ${\mathcal S}(\theta- \theta^\star)$
normalised by the posterior measure of $B_R$, is bounded on $\mathcal
{A}_0\cap\mathcal{A}_1$ by
\[
\frac{1-\Delta_{\pi}}{1+ \Delta_{\pi}}\frac{ \tilde{\mu}(
dv) } {\bar\mu(B_R)} \leq \frac{d p({\mathcal S}(\theta- \theta^\star) \mid Y) }{p(B_R\mid
Y)} \leq
\frac{ \bar{\mu}( dv) } {\tilde\mu(B_R)} \frac
{1+\Delta_{\pi}}{1- \Delta_{\pi}}.
\]
Therefore, the first term in (\ref{eqTriangleTV}) is bounded on
$\mathcal{A}_0\cap\mathcal{A}_1$ by
\begin{eqnarray*}
&&\bigl\| \mathbb{P}_{{\mathcal S}(\theta- \theta^\star) \mid Y} \mathbf {1}_{B_R}- \mu^\star
\mathbf{1}_{ B_R}\bigr\|_{\TV}
\\
&&\qquad \leq 2\int_{B_R}
\biggl[ \frac{\mathbb{P}(dv \mid Y) \mu^\star(B_R)}{
\mathbb{P}(B_R\mid Y) \mu^\star(dv)} -1 \biggr]_+ \frac{\mu^\star
(dv) }{\mu^\star(B_R)}
\\
&&\qquad \leq 2 \int_{B_R} \biggl[ \frac{ \bar\mu( dv) } {\tilde\mu
(B_R)}
\frac{\mu^\star(B_R)} {\mu^\star(dv)}\frac{(1+\Delta_{\pi
})}{(1- \Delta_{\pi})} -1 \biggr]_+ \frac{ \mu^\star(dv) }{\mu
^\star(B_R)}.
\end{eqnarray*}
Define $\mu_0 (dv) = \mu(dv; a_1, \alpha, \Omega_{00} a_0, \Omega
_{00})$. Then
\[
\frac{\mu^\star(dv)}{\mu^\star(B_R)} = \frac{\mu_0(dv)}{\mu
_0(B_R)}
\]
and
\[
\frac{\bar\mu( dv)}{\mu_0( dv)}=
\exp\bigl\{ \delta_{*1} \mathbf{1}^T v_1 +
\delta_{*0} \| v_{0}\|^2/2 \bigr\},
\]
which implies
\begin{eqnarray*}
&& \frac{ \bar\mu( dv) }{ \mu_0( dv)} \frac{\mu_0(B_R)}{\tilde\mu(B_R)}
\\
&&\qquad = \exp\bigl\{ \delta_{*0}
\|v_{0} \|^2/2 + \delta_{*1}
\mathbf{1}^T v_1 \bigr\}
\\
&&\quad\qquad{}\times
\Biggl({\int_{B_R} \prod_{i\in T_0^\star\cup T_1}
v_{i}^{\alpha_i-1} \exp\bigl\{- a_1^T v_1 \bigr\} \exp\bigl\{ - \bigl\|\Omega_{00}^{
1/2} v_{0} \bigr\|^2/2 + v_{0}^T \Omega_{00} a_0 \bigr\} \,dv}\Biggr)
\\
&&\hspace*{48pt}\bigg/\Biggl(\int_{B_R}\prod_{i\in T_0^\star\cup T_1} v_{i}^{\alpha_i-1} \exp\bigl\{- (a_1+\delta
_{*1} \mathbf{1})^T v_1 - \bigl\|\widetilde\Omega_{00}^{1/2} v_{0} \bigr\|^2/2
\\[-7pt]
&&\hspace*{259pt}{}
+ v_{0}^T \Omega_{00} a_0 \bigr\} \,dv \Biggr).
\end{eqnarray*}
To show that this expression is greater than 1, it is sufficient to
show that for any $\mathcal{B}\subseteq\{ v_{0}\dvtx  (v_{0}^T,
v_1^T)^T\in B_R\}$, the following expression is positive:
\begin{eqnarray*}
&& \int_{\mathcal{B}} \prod_{i\in T_0^\star}
w_{i}^{\alpha_i-1} e^{w^T \Omega_{00}a_0 - \|\Omega_{00}^{1/2} w \|^2/2} \,d w -\int_{\mathcal{B}}
\prod_{i\in T_0^\star} w_{i}^{\alpha_i-1}
e^{w^T
\Omega_{00}a_0 - \|\widetilde\Omega_{00}^{1/2} w \|^2/2 } \,dw
\\
&&\qquad = \int_{\mathcal{B}} \prod_{i\in T_0^\star}
w_{i}^{\alpha_i-1} e^{
- \|\widetilde\Omega_{00}^{1/2} w \|^2/2 + w^T \Omega_{00}a_0}\bigl[ \exp \bigl\{
\delta_{*0}\| w \|^2/2 \bigr\} -1\bigr] \,dw >0
\end{eqnarray*}
which is the case.
Thus, on $\mathcal{A}_0\cap\mathcal{A}_1$, $( \bar\mu( dv) /
\mu_0( dv)) (\mu_0(B_R)/\tilde\mu(B_R)) \geq1$ and hence
\begin{eqnarray*}
&& \bigl\| \mathbb{P}_{{\mathcal S}(\theta- \theta^\star) \mid Y} \mathbf {1}_{B_R}- \mu^\star
\mathbf{1}_{ B_R}\bigr\|_{\TV}
\\
&&\qquad \leq 2 \int_{B_R}
\biggl[ \frac{ \bar\mu( dv) } {\tilde\mu(B_R)} \frac{\mu^\star(B_R)} {\mu^\star(dv)}\frac{(1+\Delta_{\pi
})}{(1- \Delta_{\pi})} -1 \biggr]
\frac{ \mu^\star(dv) }{\mu^\star(B_R)}
\\
&&\qquad = 2 \biggl[ \frac{ \bar\mu(B_R) } {\tilde\mu(B_R)}\frac
{(1+\Delta_{\pi})}{(1- \Delta_{\pi})} -1 \biggr]
\\
&&\qquad = 2
\frac{ \bar\mu(B_R) -\tilde\mu(B_R) } {\tilde\mu
(B_R)}\frac{(1+\Delta_{\pi})}{(1- \Delta_{\pi})} + 2 \biggl[ \frac
{(1+\Delta_{\pi})}{(1- \Delta_{\pi})} -1 \biggr].
\end{eqnarray*}
The difference of measures $ \bar\mu(B_R) -\tilde\mu(B_R)$ is
bounded by
\begin{eqnarray*}
&& \int_{B_R}\prod_{i\in T_0^\star\cup T_1}
v_{i}^{\alpha_i-1} e^{
- v_{0}^T\widetilde\Omega_{00} v_{0} /2 + v_0^T \Omega_{00} a_0 -
\tilde{a} v_1 } \bigl[ e^{ \delta_{*0} \|v_0\|^2/2 + \delta
_{*1} 1_{p_1}^T v_1 } -1
\bigr] \,dv
\\
&&\qquad \leq \int_{B_R} \prod_{i \in T_0^\star\cup T_1}
v_{i}^{\alpha
_i-1} \bigl[\delta_{*0} \|v_0
\|^2/2 + \delta_{*1} 1_{p_1}^T
v_1\bigr] e^{ -
v_{0}^T\widebar\Omega_{00} v_{0} /2 + v_0^T \Omega_{00} a_0 - \bar{a}
v_1 } \,dv
\\
&&\qquad \leq \Biggl[\delta_{*0} E_{\Phi} + \delta_{*1}
\sum_{j=1}^{p_1} (\alpha_{1,j}/
\bar{a}_j ) \Biggr] \bar{\mu}\bigl(\mathcal{V}^\star\bigr)
\end{eqnarray*}
due to the inequality $e^x-1\leq x e^x$ for $x>0$, where $E_{\Phi}$ is
defined by
%
\begin{equation}
\label{eqdefEPhi} E_{\Phi} = 0.5 \int_{\mathcal{V}_0} \|w
\|^2 \mathcal{PTN}_{p_0}\bigl(dw; \widebar\Omega_{00}^{-1}
\Omega_{00} a_0, \widebar\Omega_{00}^{-1},
p_0^\star,\alpha_0\bigr),
\end{equation}
which is finite. Therefore,
\begin{eqnarray*}
&& \bigl\| \mathbb{P}_{{\mathcal S}(\theta- \theta^\star) \mid Y} \mathbf {1}_{B_R}- \mu^\star
\mathbf{1}_{ B_R}\bigr\|_{\TV}
\\
&&\qquad \leq\frac{2\bar{\mu
}(\mathcal{V}^\star)}{\tilde\mu(B_R)}
\frac{(1+\Delta_{\pi})}{(1- \Delta_{\pi})} \Biggl[\delta_{*0} E_{\Phi} +
\delta_{*1} \sum_{j=1}^{p_1}
\frac{\alpha_{1,j}}{\bar{a}_j } \Biggr]
+ \frac{4\Delta_{\pi}}{1- \Delta_{\pi}},
\end{eqnarray*}
which goes to zero since $\delta_{*k} \to0$ and $\Delta_{\pi} \to
0$ as $\sigma\to0$.
For small $\sigma$ and hence large $R_0$ and $R_1$, the ratios $\bar
{\mu}(\mathcal{V}^\star)/\tilde\mu(\mathcal{V}^\star)$ and
\[
\frac{\tilde\mu(B_R)}{\tilde\mu(\mathcal{V}^\star)} = \mathcal{PTN}_{p_0}\bigl( B_{2}(0,
R_0); \widetilde\Omega_{00}^{-1}\Omega
_{00} a_0, \widetilde\Omega_{00}^{-1},
p^\star_0, \alpha_0\bigr) \prod
_{j=1}^{p_1} \Gamma\bigl((0,R_1);
\alpha_{1,j}, \tilde{a}_{j}\bigr)
\]
are close to 1. Therefore, $\| \mathbb{P}_{{\mathcal S}(\theta-
\theta^\star) \mid Y}\mathbf{1}_{B_R} - \mu^\star\mathbf{1}_{B_R}
\|_{\TV} \to0$ as $\sigma\to0$.

The second term in (\ref{eqTriangleTV}) is bounded by
$\|\mu^\star- \mu^\star\mathbf{1}_{B_R}\|_{\TV} \leq2 \mu^\star
(\overline{B_R}) \to0$
as $R_0, R_1\to\infty$, since the set $B_R$ converges to $\mathcal
{V}^\star$ by Lemma~\ref{lemLimitX}.

The third term in (\ref{eqTriangleTV}) is bounded by
\[
\bigl\| \mathbb{P}_{({\mathcal S}(\theta- \theta^\star) \mid Y)} \mathbf {1}_{B_R}- \mathbb{P}_{{\mathcal S}(\theta- \theta^\star) \mid Y}
\bigr\| _{\TV} \leq2 \mathbb{P}_{{\mathcal S}(\theta- \theta^\star) \mid
Y}(\overline{B_R})
\leq\frac{2\Delta_0(\delta)}{ C_{\pi} (1-
\Delta_{\pi}) \tilde{\mu}(B_R)},
\]
where $\Delta_0( \delta)$ is defined by (\ref{eqDefD0}).
By Assumption~\ref{assL}, with probability $\to1$, $\Delta_0(\delta) \to0$
as $\sigma\to0$; also, $\tilde{\mu}(B_R) \to\mu_0(\mathcal
{V}^\star) >0$.

Combining these bounds, we have that on $\mathcal{A}_0\cap\mathcal{A}_1$,
\begin{eqnarray*}
&&\bigl\| \mathbb{P}_{{\mathcal S}(\theta- \theta^\star) \mid Y} - \mu ^\star\bigr\|_{\TV}
\\
&&\qquad \leq 2
\mu^\star(\overline{B_R}) +2 \bigl[C_{\pi} (1-
\Delta_{\pi}) \tilde{\mu}(B_R)\bigr]^{-1}
\Delta_0(\delta)
\\
&&\quad\qquad{}+ 2\frac{\bar{\mu}(\mathcal{V}^\star)}{\tilde\mu(B_R)} \frac{(1+\Delta_{\pi})}{(1- \Delta_{\pi})} \Biggl[\delta_{*0}
E_{\Phi} + \delta_{*1} \sum_{j=1}^{p_1}
(\alpha_{1,j}/\bar{a}_j ) \Biggr] + \frac{4 \Delta_{\pi} }{(1- \Delta_{\pi})}
\to0 
\end{eqnarray*}
and \mbox{$\mathbb{P}_{\theta^\star, \sigma}(\mathcal{A}_0\cap\mathcal
{A}_1) \to1$} as $\sigma\to0$ due to Assumption~\ref{assS}, which gives the
statement of the theorem.
\end{pf*}
%

\begin{pf*}{Proof of Proposition~\ref{propBvMnonasympt}}
In the proof of Theorem~\ref{thPostApproxBvM}, we derived the
following upper bound on event $\mathcal{A}$:
\[
\bigl\| \mathbb{P}_{{\mathcal S}(\theta- \theta^\star) \mid Y} - \mu ^\star\bigr\|_{\TV} \leq 2
\mu^\star(\overline{B_R}) + C_{\Delta}
\Delta_0(\delta) + 2C_{0}\delta_{*0} +
2C_{1}\delta_{*1} + C_{2}
\Delta_{\pi},
\]
where $C_{\Delta} = 2[C_{\pi} (1 -\Delta_{\pi})\tilde{\mu
}(B_R)]^{-1}$, $C_{2} = 4/(1- \Delta_{\pi})$, $C_{0} = C_A E_{\Phi}$
with $E_\Phi$ defined by (\ref{eqdefEPhi}), $C_{1} = C_A \sum_{j=1}^{p_1} \alpha_{1,j}/(a_{1,j}- \delta_{*1})$
and with $B_{R,0}=\break B_{2, p_0}(0, R_0)\cap\mathcal{V}_0$,
\[
C_{A} = \frac{\bar\mu(\mathcal{V}^\star)}{\tilde\mu(B_R)} \frac{(1+\Delta_{\pi})}{(1- \Delta_{\pi})} = \frac{\bar\mu_{p_0}(\mathcal{V}_0)}{\tilde\mu_{p_0}(B_{R,0})}
\prod_{j=1}^{p_1} \biggl[ \frac{a_{1,j} + \delta_{*1}}{a_{1,j} -
\delta_{*1}}
\biggr]^{\alpha_{1,j}} \frac{(1+\Delta_{\pi})}{(1- \Delta_{\pi})},
\]
where $\mu_{p_0}(\mathcal{B}_0) = \int_{\mathcal{B}_0 \times
[0,\infty)^{p_1}} \mu(dv)$ for a measure $\mu$, $\mathcal{B}_0
\subset\mathcal{V}_0$.
If $S_0^\star=\varnothing$,
\begin{eqnarray*}
E_{\Phi} &=& \bigl\|\widebar\Omega_{00}^{-1}
\Omega_{00} a_0\bigr\|^2/2 + \operatorname{trace}
\bigl(\widebar\Omega_{00}^{-1}\bigr)/2,
\\
\frac{\bar\mu_{p_0}(\mathcal{V}_0)}{\tilde\mu_{p_0}(B_{R,0})} &=& \frac{ e^{ \delta_{*0} a_0^T \Omega_{00} \widebar\Omega
_{00}^{-1}\widetilde\Omega_{00}^{-1}\Omega_{00} a_0} [\det(\widebar
\Omega_{00}^{-1}\widetilde\Omega_{00})]^{1/2}}{\mathcal
{TN}(B_{R,0}; \widetilde\Omega_{00}^{-1}\Omega_{00}a_0, \widetilde
\Omega_{00}^{-1})}.
\end{eqnarray*}
We bound the term $\mu^\star (\overline{B_R} )$ by
\begin{eqnarray*}
\mu^\star(\overline{B_R}) &=& 1 - \mu^\star_{p_0}(B_{R,0})
\prod_{j=1}^{p_1} \Gamma \biggl( \biggl(0,
\frac{\delta_1}{\sigma^2} \biggr); \alpha_{1,j}, a_{1,j} \biggr)
\\
&\leq& \mu^\star_{p_0}(\overline{B_{R,0}}) +1 -
\prod_{j=1}^{p_1} \Gamma \biggl( \biggl(0,
\frac{\delta_1}{\sigma^2} \biggr); \alpha _{1,j}, a_{1,j} \biggr)
\end{eqnarray*}
using the inequality $1-xy \leq1- x + 1-y$ for $x,y \in(0,1)$. We can
also use
\begin{eqnarray*}
1 -\prod_{j=1}^{p_1} \Gamma \biggl(
\biggl(0,\frac{\delta_1}{\sigma
^2} \biggr); \alpha_{1,j}, a_{1,j}
\biggr) &\leq& p_1\Bigl[1 - \min_{j}\Gamma
\bigl(\bigl(0,\delta_1/\sigma^2\bigr);
\alpha_{1,j}, a_{1,j} \bigr)\Bigr]
\\
&=& p_1 \max_{j}\Gamma \bigl(\bigl(
\delta_1/\sigma^2,\infty\bigr); \alpha_{1,j},
a_{1,j} \bigr),
\end{eqnarray*}
and, changing to polar coordinates and denoting $p_{\alpha0} = p_0 +
\sum_{j\in T_0^\star}(\alpha_{0,j}-1)$ and $W =\{ w\in\mathbb
{R}^{p_0}\dvtx  \|w\|_2^2 =1, w_{j}>0$ for $j\in T_0^\star
\}$, we have
\begin{eqnarray*}
\mu^\star_{p_0} (\overline{B_{R,0}} ) &\leq& \mu
_0\bigl(\mathcal{V}^\star\bigr)\int_{R_0}^{\infty}
r^{p_{\alpha0}-1} e^{-
\lambda_{\min}(\Omega_{00}) (r- \|a_0\|)^2/2}\,dr \int_{W} \prod
_{j\in T_0^\star} w_j^{\alpha_{0,j}-1} \,dw
\\
&\leq& C_{\alpha0} \Gamma \bigl(\bigl(\bigl(\delta_0/\sigma-
\|a_0\| \bigr)^2/2,\infty\bigr); p_{\alpha0}/2,
\lambda_{\min}(\Omega_{00}) \bigr)
\end{eqnarray*}
under the assumption that $R_0=\delta_0/\sigma> \|a_0(\omega)\|$ where
\[
C_{\alpha0} = \mu_0\bigl(\mathcal{V}^\star\bigr)
2^{-p_0^\star+1.5p_{\alpha
0}} \bigl[ \lambda_{\min}(\Omega_{00})
\bigr]^{p_{\alpha0}/2}\pi ^{(p_0-p^\star_0)/2}\prod_{i\in T_0^\star}
\Gamma(\alpha_{0,i}/2).
\]

Collecting conditions on $\delta_k$ used in the proof of Theorem~\ref
{thPostApproxBvM}, we have conditions~(\ref{eqNonAsympCond}). Thus,
we have the required inequality on the event $\mathcal{A}$.
\end{pf*}


\subsection{Auxiliary results}\label{secProofAux}
\mbox{}

%
\begin{pf*}{Proof of Lemma~\ref{lemLimitX}}
Due to Assumption~\ref{assB} and the fact that $\theta^\star_{S^\star_0 \cup
S_1} =0$, the set $D^{-1} U(\Theta^\star(\delta) - \theta^\star)$ contains
\begin{eqnarray*}
&& B_{2, p_0} \biggl(0, \frac{\delta_0}{\sigma} \biggr)\times B_{\infty, p_1}
\biggl(0, \frac{\delta_1}{\sigma^2} \biggr) \cap \biggl(- \frac{c_0}{\sigma},
\frac{c_0}{\sigma} \biggr)^{p_0-p_0^\star} \times \biggl[0, \frac{c_0}{\sigma}
\biggr)^{p_0^\star} \times \biggl[0, \frac{c_1}{\sigma^2} \biggr)^{p_1}
\\
&&\qquad = \bigl\{v\dvtx  v\in B_{2, p_0} (0, \delta_0/\sigma)\mbox{ and }
v_{T_0^\star} \geq0\bigr\} \times\bigl[0, \delta_1/
\sigma^2 \bigr)^{p_1},
\end{eqnarray*}
where $T_0^\star=\{p_0-p_0^\star+1,\ldots,p_0\}$. These sets
monotonically increase to $\mathcal{V}^\star= \mathbb
{R}^{p_0-p^\star_0}\times\mathbb{R}^{p^\star_0+p_1}_+$ as $\sigma
\to0$ due to the assumption $\delta_0/\sigma\to\infty$ and $\delta
_1/\sigma^2 \to\infty$; this implies the statement of the lemma.
\end{pf*}


\begin{pf*}{Proof of Lemma~\ref{lemAssumptionL}}
Under the assumptions of the lemma, for small enough~$\sigma$, with
$\tilde{\delta}_0 = \delta_0/\sqrt{p_0}$, we have that
\begin{eqnarray*}
&& \frac{1} { C_{\pi0}(\delta)} \int_{\Theta\setminus\Theta^\star
(\delta)} e^{ (\ell_{y }(\theta) - \ell_{y }(\theta^\star
))/\sigma^2} \pi(d\theta)
\\
&&\qquad \leq \sum_{j \in S_0} \int_{\tilde{\delta}_0}^{\infty}
e^{-C_{\delta0}
v_j/\sigma^2} \,dv_j
\\
&&\quad\qquad{}+ \sum_{j\in S_0\setminus S_0^\star} \int_{0}^{\theta^\star_j -
\tilde{\delta}_0}
\theta_j^{\bolds\alpha_j-1} e^{ -C_{\delta0}
|\theta
_j-\theta^\star_j|/\sigma^2} \,d\theta_j
\\
&&\quad\qquad{}+
\sum_{j \in S_1} \int_{\delta_1}^{\infty}
e^{-C_{\delta1}
v_j/\sigma^2} \,dv_j
\\
&&\qquad \leq\sum_{j\in S_0\setminus S_0^\star} \sigma^{\bolds\alpha_j}
e^{-C_{\delta0}(\theta^\star_j - \sigma)/\sigma^2} + \frac{p_0\sigma^2}{C_{\delta0}} e^{-C_{\delta} \delta_0/[\sqrt
{p_0}\sigma^2]} + \frac{p_1\sigma^2}{C_{\delta1}}
e^{-C_{\delta1}
\delta_1/\sigma^2}
\\
&&\quad\qquad{}+ \sum_{j\in S_0\setminus S_0^\star} \bigl[\sigma^{\bolds\alpha
_j-1}I(\bolds
\alpha_j<1) + {\theta^\star_j}^{\bolds\alpha
_j-1}I(
\bolds\alpha _j\geq1) \bigr] \frac{\sigma^2}{C_{\delta0}} e^{-C_{\delta
0}\tilde{\delta}_0/ \sigma^2}
\\
&&\qquad \leq C \bigl[ \sigma^{\min_{j}(\bolds\alpha_j)} + \sigma\bigr] e^{-C_{\delta0}
\delta_0/[\sqrt{p_0}\sigma^2]} +
p_1 e^{-C_{\delta1} \delta
_1/\sigma^2} \sigma^2/C_{\delta1}
\end{eqnarray*}
for a constant $C$.
This implies that, with $J_\sigma= \sigma^{-\sum_{j\in S_0} \bolds
\alpha
_{j} - 2\sum_{j\in S_1}\bolds\alpha_{j}}$,
\[
\Delta_0 (\delta) \leq C_{\pi0}(\delta)J_\sigma
\biggl[ C\bigl[ \sigma^{\min_{j}(\bolds\alpha
_j)} + \sigma\bigr] e^{-C_{\delta0} \delta_0/[\sqrt{p_0}\sigma^2]} +
\frac{
p_1 \sigma^2}{C_{\delta1}} e^{-C_{\delta1} \delta_1/\sigma^2} \biggr] \to0
\]
as $\sigma\to0$ under the assumptions of the lemma.
\end{pf*}
\end{appendix}



%

\printaddresses

\begin{thebibliography}{28}
\bibitem[\protect\citeauthoryear{Barron, Schervish and
Wasserman}{1999}]{Barron99}
%
\begin{barticle}[mr]
\bauthor{\bsnm{Barron},~\bfnm{Andrew}\binits{A.}},
\bauthor{\bsnm{Schervish},~\bfnm{Mark~J.}\binits{M.~J.}} \AND
\bauthor{\bsnm{Wasserman},~\bfnm{Larry}\binits{L.}}
(\byear{1999}).
\btitle{The consistency of posterior distributions in nonparametric problems}.
\bjournal{Ann. Statist.}
\bvolume{27}
\bpages{536--561}.
\bid{doi={10.1214/aos/1018031206}, issn={0090-5364}, mr={1714718}}
\end{barticle}
%
\bptok{imsref}%
\endbibitem

\bibitem[\protect\citeauthoryear{Bertsekas}{2003}]{Bertsekas}
%
\begin{bbook}[author]
\bauthor{\bsnm{Bertsekas},~\bfnm{D.~P.}\binits{D.~P.}}
(\byear{2003}).
\btitle{Convex Analysis and Optimization}.
\bpublisher{Athena Scientific and Tsinghua Univ. Press},
\blocation{Belmont, MA}.
\end{bbook}
%
\bptok{imsref}%
\endbibitem

\bibitem[\protect\citeauthoryear{Besag}{1986}]{Besag}
%
\begin{barticle}[mr]
\bauthor{\bsnm{Besag},~\bfnm{Julian}\binits{J.}}
(\byear{1986}).
\btitle{On the statistical analysis of dirty pictures}.
\bjournal{J. Roy. Statist. Soc. Ser. B}
\bvolume{48}
\bpages{259--302}.
\bid{issn={0035-9246}, mr={0876840}}
\bptnote{check related}%
\end{barticle}
%
\bptok{imsref}%
\endbibitem

\bibitem[\protect\citeauthoryear{Bochkina}{2013}]{BochkinaIP}
%
\begin{barticle}[mr]
\bauthor{\bsnm{Bochkina},~\bfnm{Natalia}\binits{N.}}
(\byear{2013}).
\btitle{Consistency of the posterior distribution in generalized linear
inverse problems}.
\bjournal{Inverse Problems}
\bvolume{29}
\bpages{095010, 43}.
\bid{doi={10.1088/0266-5611/29/9/095010}, issn={0266-5611}, mr={3094485}}
\end{barticle}
%
\bptok{imsref}%
\endbibitem

\bibitem[\protect\citeauthoryear{Chernozhukov and Hong}{2004}]{ChernozhukovHong}
%
\begin{barticle}[mr]
\bauthor{\bsnm{Chernozhukov},~\bfnm{Victor}\binits{V.}} \AND
\bauthor{\bsnm{Hong},~\bfnm{Han}\binits{H.}}
(\byear{2004}).
\btitle{Likelihood estimation and inference in a class of nonregular
econometric models}.
\bjournal{Econometrica}
\bvolume{72}
\bpages{1445--1480}.
\bid{doi={10.1111/j.1468-0262.2004.00540.x}, issn={0012-9682}, mr={2077489}}
\end{barticle}
%
\bptok{imsref}%
\endbibitem

\bibitem[\protect\citeauthoryear{Douc et~al.}{2011}]{DoucHMM}
%
\begin{barticle}[mr]
\bauthor{\bsnm{Douc},~\bfnm{Randal}\binits{R.}},
\bauthor{\bsnm{Moulines},~\bfnm{Eric}\binits{E.}},
\bauthor{\bsnm{Olsson},~\bfnm{Jimmy}\binits{J.}} \AND
\bauthor{\bparticle{van} \bsnm{Handel},~\bfnm{Ramon}\binits{R.}}
(\byear{2011}).
\btitle{Consistency of the maximum likelihood estimator for general
hidden {M}arkov models}.
\bjournal{Ann. Statist.}
\bvolume{39}
\bpages{474--513}.
\bid{doi={10.1214/10-AOS834}, issn={0090-5364}, mr={2797854}}
\end{barticle}
%
\bptok{imsref}%
\endbibitem

\bibitem[\protect\citeauthoryear{Dudley and Haughton}{2002}]{HDudley}
%
\begin{barticle}[mr]
\bauthor{\bsnm{Dudley},~\bfnm{R.~M.}\binits{R.~M.}} \AND
\bauthor{\bsnm{Haughton},~\bfnm{D.}\binits{D.}}
(\byear{2002}).
\btitle{Asymptotic normality with small relative errors of posterior
probabilities of half-spaces}.
\bjournal{Ann. Statist.}
\bvolume{30}
\bpages{1311--1344}.
\bid{doi={10.1214/aos/1035844978}, issn={0090-5364}, mr={1936321}}
\end{barticle}
%
\bptok{imsref}%
\endbibitem

\bibitem[\protect\citeauthoryear{Erkanli}{1994}]{Erkanli}
%
\begin{barticle}[mr]
\bauthor{\bsnm{Erkanli},~\bfnm{Alaattin}\binits{A.}}
(\byear{1994}).
\btitle{Laplace approximations for posterior expectations when the mode
occurs at the boundary of the parameter space}.
\bjournal{J. Amer. Statist. Assoc.}
\bvolume{89}
\bpages{250--258}.
\bid{issn={0162-1459}, mr={1266297}}
\end{barticle}
%
\bptok{imsref}%
\endbibitem

\bibitem[\protect\citeauthoryear{Geman and Geman}{1984}]{GemanG}
%
\begin{barticle}[pbm]
\bauthor{\bsnm{Geman},~\bfnm{S.}\binits{S.}} \AND
\bauthor{\bsnm{Geman},~\bfnm{D.}\binits{D.}}
(\byear{1984}).
\btitle{Stochastic relaxation, gibbs distributions, and the Bayesian
restoration of images}.
\bjournal{IEEE Trans. Pattern Anal. Mach. Intell.}
\bvolume{6}
\bpages{721--741}.
\bid{issn={0162-8828}, pmid={22499653}}
\end{barticle}
%
\bptok{imsref}%
\endbibitem

\bibitem[\protect\citeauthoryear{Ghosal, Ghosh and Samanta}{1995}]{GGS95}
%
\begin{barticle}[mr]
\bauthor{\bsnm{Ghosal},~\bfnm{Subhashis}\binits{S.}},
\bauthor{\bsnm{Ghosh},~\bfnm{Jayanta~K.}\binits{J.~K.}} \AND
\bauthor{\bsnm{Samanta},~\bfnm{Tapas}\binits{T.}}
(\byear{1995}).
\btitle{On convergence of posterior distributions}.
\bjournal{Ann. Statist.}
\bvolume{23}
\bpages{2145--2152}.
\bid{doi={10.1214/aos/1034713651}, issn={0090-5364}, mr={1389869}}
\end{barticle}
%
\bptok{imsref}%
\endbibitem

\bibitem[\protect\citeauthoryear{Ghosal and Samanta}{1995}]{GS95}
%
\begin{barticle}[mr]
\bauthor{\bsnm{Ghosal},~\bfnm{Subhashis}\binits{S.}} \AND
\bauthor{\bsnm{Samanta},~\bfnm{Tapas}\binits{T.}}
(\byear{1995}).
\btitle{Asymptotic behaviour of {B}ayes estimates and posterior
distributions in multiparameter nonregular cases}.
\bjournal{Math. Methods Statist.}
\bvolume{4}
\bpages{361--388}.
\bid{issn={1066-5307}, mr={1372011}}
\end{barticle}
%
\bptok{imsref}%
\endbibitem

\bibitem[\protect\citeauthoryear{Ghosh, Ghosal and Samanta}{1994}]{GGS94}
%
\begin{bincollection}[mr]
\bauthor{\bsnm{Ghosh},~\bfnm{Jayanta~K.}\binits{J.~K.}},
\bauthor{\bsnm{Ghosal},~\bfnm{Subhashis}\binits{S.}} \AND
\bauthor{\bsnm{Samanta},~\bfnm{Tapas}\binits{T.}}
(\byear{1994}).
\btitle{Stability and convergence of the posterior in non-regular problems}.
In \bbooktitle{Statistical Decision Theory and Related Topics, {V}
({W}est {L}afayette, {IN},~1992)}
(\beditor{\bfnm{S.~S.}\binits{S.~S.}~\bsnm{Gupta}} \AND
  \beditor{\bfnm{J.~O.}\binits{J.~O.}~\bsnm{Berger}}, eds.)
\bpages{183--199}.
\bpublisher{Springer},
\blocation{New York}.
\bid{mr={1286304}}
\end{bincollection}
%
\bptok{imsref}%
\endbibitem

\bibitem[\protect\citeauthoryear{Green}{1990}]{Green}
%
\begin{barticle}[author]
\bauthor{\bsnm{Green},~\bfnm{Peter~J.}\binits{P.~J.}}
(\byear{1990}).
\btitle{Bayesian reconstructions from emission tomography data using a
modified EM algorithm}.
\bjournal{IEEE Trans. Med. Imag.}
\bvolume{9}
\bpages{84--93}.
\end{barticle}
%
\bptok{imsref}%
\endbibitem

\bibitem[\protect\citeauthoryear{Hirano and Porter}{2003}]{HiranoPorter}
%
\begin{barticle}[mr]
\bauthor{\bsnm{Hirano},~\bfnm{Keisuke}\binits{K.}} \AND
\bauthor{\bsnm{Porter},~\bfnm{Jack~R.}\binits{J.~R.}}
(\byear{2003}).
\btitle{Asymptotic efficiency in parametric structural models with
parameter-dependent support}.
\bjournal{Econometrica}
\bvolume{71}
\bpages{1307--1338}.
\bid{doi={10.1111/1468-0262.00451}, issn={0012-9682}, mr={2000249}}
\end{barticle}
%
\bptok{imsref}%
\endbibitem

\bibitem[\protect\citeauthoryear{Ibragimov and Has'minski{\u\i}}{1981}]{IH81}
%
\begin{bbook}[mr]
\bauthor{\bsnm{Ibragimov},~\bfnm{I.~A.}\binits{I.~A.}} \AND
\bauthor{\bsnm{Has'minski{\u\i}},~\bfnm{R.~Z.}\binits{R.~Z.}}
(\byear{1981}).
\btitle{Statistical Estimation: Asymptotic Theory}.
\bpublisher{Springer},
\blocation{New York}.
\bid{mr={0620321}}
\end{bbook}
%
\bptok{imsref}%
\endbibitem

\bibitem[\protect\citeauthoryear{Johnstone and Silverman}{1990}]{JSPET}
%
\begin{barticle}[mr]
\bauthor{\bsnm{Johnstone},~\bfnm{Iain~M.}\binits{I.~M.}} \AND
\bauthor{\bsnm{Silverman},~\bfnm{Bernard~W.}\binits{B.~W.}}
(\byear{1990}).
\btitle{Speed of estimation in positron emission tomography and related
inverse problems}.
\bjournal{Ann. Statist.}
\bvolume{18}
\bpages{251--280}.
\bid{doi={10.1214/aos/1176347500}, issn={0090-5364}, mr={1041393}}
\end{barticle}
%
\bptok{imsref}%
\endbibitem

\bibitem[\protect\citeauthoryear{LeCam}{1953}]{LeCam}
%
\begin{barticle}[mr]
\bauthor{\bsnm{LeCam},~\bfnm{Lucien}\binits{L.}}
(\byear{1953}).
\btitle{On some asymptotic properties of maximum likelihood estimates
and related {B}ayes' estimates}.
\bjournal{Univ. California Publ. Statist.}
\bvolume{1}
\bpages{277--329}.
\bid{mr={0054913}}
\end{barticle}
%
\bptok{imsref}%
\endbibitem

\bibitem[\protect\citeauthoryear{Le~Cam and Yang}{1990}]{LeCamYang}
%
\begin{bbook}[mr]
\bauthor{\bsnm{Le Cam},~\bfnm{Lucien}\binits{L.}} \AND
\bauthor{\bsnm{Yang},~\bfnm{Grace~Lo}\binits{G.~L.}}
(\byear{1990}).
\btitle{Asymptotics in Statistics: Some Basic Concepts}.
\bpublisher{Springer},
\blocation{New York}.
\bid{doi={10.1007/978-1-4684-0377-0}, mr={1066869}}
\end{bbook}
%
\bptok{imsref}%
\endbibitem

\bibitem[\protect\citeauthoryear{Moran}{1971}]{Moran}
%
\begin{barticle}[mr]
\bauthor{\bsnm{Moran},~\bfnm{P.~A.~P.}\binits{P.~A.~P.}}
(\byear{1971}).
\btitle{Maximum-likelihood estimation in non-standard conditions}.
\bjournal{Proc. Cambridge Philos. Soc.}
\bvolume{70}
\bpages{441--450}.
\bid{mr={0290493}}
\end{barticle}
%
\bptok{imsref}%
\endbibitem

\bibitem[\protect\citeauthoryear{Nelder and Wedderburn}{1972}]{NelderW}
%
\begin{barticle}[author]
\bauthor{\bsnm{Nelder},~\bfnm{J.~A.}\binits{J.~A.}} \AND
\bauthor{\bsnm{Wedderburn},~\bfnm{R.~W.~M.}\binits{R.~W.~M.}}
(\byear{1972}).
\btitle{Generalized linear models}.
\bjournal{J. R. Stat. Soc. A, General}
\bvolume{135}
\bpages{370--384}.
\end{barticle}
%
\bptok{imsref}%
\endbibitem

\bibitem[\protect\citeauthoryear{Petrone, Rousseau and
Scricciolo}{2012}]{PetroneEB}
%
\begin{barticle}[author]
\bauthor{\bsnm{Petrone},~\bfnm{S.}\binits{S.}},
\bauthor{\bsnm{Rousseau},~\bfnm{J.}\binits{J.}} \AND
\bauthor{\bsnm{Scricciolo},~\bfnm{C.}\binits{C.}}
(\byear{2012}).
\btitle{Bayes and empirical {B}ayes: Do they merge?}
\bjournal{Biometrika}
\bvolume{99}
\bpages{1--21}.
\end{barticle}
%
\bptok{imsref}%
\endbibitem

\bibitem[\protect\citeauthoryear{Polyak}{1983}]{Polyak}
%
\begin{bbook}[mr]
\bauthor{\bsnm{Polyak},~\bfnm{B.~T.}\binits{B.~T.}}
(\byear{1983}).
\btitle{Introduction to Optimization (Vvedenie v Optimizatsiyu, in Russian)}.
\bpublisher{Nauka},
\blocation{Moscow}.
\bid{mr={0719196}}
\end{bbook}
%
\bptok{imsref}%
\endbibitem

\bibitem[\protect\citeauthoryear{Self and Liang}{1987}]{SelfLiang}
%
\begin{barticle}[mr]
\bauthor{\bsnm{Self},~\bfnm{Steven~G.}\binits{S.~G.}} \AND
\bauthor{\bsnm{Liang},~\bfnm{Kung-Yee}\binits{K.-Y.}}
(\byear{1987}).
\btitle{Asymptotic properties of maximum likelihood estimators and
likelihood ratio tests under nonstandard conditions}.
\bjournal{J. Amer. Statist. Assoc.}
\bvolume{82}
\bpages{605--610}.
\bid{issn={0162-1459}, mr={0898365}}
\end{barticle}
%
\bptok{imsref}%
\endbibitem

\bibitem[\protect\citeauthoryear{Shapiro}{2000}]{Shapiro}
%
\begin{barticle}[mr]
\bauthor{\bsnm{Shapiro},~\bfnm{Alexander}\binits{A.}}
(\byear{2000}).
\btitle{On the asymptotics of constrained local {$M$}-estimators}.
\bjournal{Ann. Statist.}
\bvolume{28}
\bpages{948--960}.
\bid{doi={10.1214/aos/1015952006}, issn={0090-5364}, mr={1792795}}
\end{barticle}
%
\bptok{imsref}%
\endbibitem

\bibitem[\protect\citeauthoryear{Spokoiny}{2012}]{Spokoiny2012}
%
\begin{barticle}[mr]
\bauthor{\bsnm{Spokoiny},~\bfnm{Vladimir}\binits{V.}}
(\byear{2012}).
\btitle{Parametric estimation. {F}inite sample theory}.
\bjournal{Ann. Statist.}
\bvolume{40}
\bpages{2877--2909}.
\bid{doi={10.1214/12-AOS1054}, issn={0090-5364}, mr={3097963}}
\end{barticle}
%
\bptok{imsref}%
\endbibitem

\bibitem[\protect\citeauthoryear{van~der Vaart}{1998}]{vdVaart}
%
\begin{bbook}[mr]
\bauthor{\bsnm{van~der Vaart},~\bfnm{A.~W.}\binits{A.~W.}}
(\byear{1998}).
\btitle{Asymptotic Statistics}.
\bpublisher{Cambridge Univ. Press},
\blocation{Cambridge}.
\bid{doi={10.1017/CBO9780511802256}, mr={1652247}}
\end{bbook}
%
\bptok{imsref}%
\endbibitem

\bibitem[\protect\citeauthoryear{Vu and Zhou}{1997}]{LRboundary}
%
\begin{barticle}[mr]
\bauthor{\bsnm{Vu},~\bfnm{H.~T.~V.}\binits{H.~T.~V.}} \AND
\bauthor{\bsnm{Zhou},~\bfnm{S.}\binits{S.}}
(\byear{1997}).
\btitle{Generalization of likelihood ratio tests under nonstandard conditions}.
\bjournal{Ann. Statist.}
\bvolume{25}
\bpages{897--916}.
\bid{doi={10.1214/aos/1031833677}, issn={0090-5364}, mr={1439327}}
\end{barticle}
%
\bptok{imsref}%
\endbibitem

\bibitem[\protect\citeauthoryear{Weir}{1997}]{Weir}
%
\begin{barticle}[author]
\bauthor{\bsnm{Weir},~\bfnm{Iain~S.}\binits{I.~S.}}
(\byear{1997}).
\btitle{Fully {B}ayesian reconstructions from single-photon emission
computed tomography data}.
\bjournal{J. Amer. Statist. Assoc.}
\bvolume{92}
\bpages{49--60}.
\end{barticle}
%
\bptok{imsref}%
\endbibitem
\end{thebibliography}
\end{document}